\documentclass[11pt]{amsart}

\usepackage{amssymb}
\usepackage{latexsym}
\usepackage{amsmath}
\usepackage{amscd}
\usepackage{array}

\def\c{{\rm C}}
\def\C{{\mathbb C}}
\def\r{{\mathcal R}}

\def\z{{\mathbb Z}}

\def\Q{{\mathcal Q}}
\def\V{{\rm V}}

\begin{document}

\newtheorem{theorem}{Theorem}[section]
\newtheorem{proposition}{Proposition}[section]
\newtheorem{definition}{Definition}[section]
\newtheorem{corollary}{Corollary}[section]
\newtheorem{lemma}{Lemma}[section]
\newtheorem{conjecture}{Conjecture}[section]
\newtheorem{question}{Question}[section]

\title[Quantum groups acting on 4 points]{Quantum groups acting on 4 points}

\author{Teodor Banica}
\address{Laboratoire de Math\'ematiques, Universit\'e Paul Sabatier, Toulouse III, 118 route de Narbonne, 31062 Toulouse, France} 
\email{Teodor.Banica@math.ups-tlse.fr}

\author{Julien Bichon}
\address{Laboratoire de Math\'ematiques, Universit\'e Blaise Pascal, Clermont-Ferrand II, 
Campus des C\'ezeaux, 63177 Aubi\`ere Cedex, France}
\email{Julien.Bichon@math.univ-bpclermont.fr}

\subjclass[2000]{20G42,16W30}
\keywords{Quantum permutation group}

\begin{abstract}
We classify the compact quantum groups acting on 4 points. These are the
quantum subgroups of the quantum permutation group $\mathcal Q_4$. Our main
tool is a new presentation for the algebra $\rm C(\mathcal Q_4)$,
corresponding to an isomorphism of type $\mathcal Q_4\simeq SO_{-1}(3)$. The 
quantum subgroups of $\mathcal Q_4$ are subject to a McKay type correspondence, that we describe at the level of algebraic invariants.
\end{abstract}

\maketitle

\section{Introduction}

After the seminal work of Drinfeld \cite{dr} and Woronowicz \cite{wo}, an influencal 
treatment in quantum group theory was Manin's book \cite{ma}.
Manin proposed to construct quantum groups as quantum symmetry
groups of quantum spaces, which  were seen there as dual objects of quadratic algebras.
Later on, Wang \cite{wa} studied quantum symmetry groups of finite quantum spaces,
corresponding to finite-dimensional $\c^*$-algebras, answering in this way 
a question of Connes on the existence of such quantum groups.
One surprising conclusion of Wang's work is that there exist infinite compact
quantum groups acting faithfully on $n$ points, provided $n\geq 4$.
The biggest such quantum group, denoted here by $\mathcal Q_n$, is called
the quantum permutation group on $n$ points. In other words, the symmetric group
$S_n$ has $\mathcal Q_n$ as an infinite quantum analogue if $n\geq 4$. 
  
Soon after Wang's discovery, the representation theory of the quantum group
$\Q_n$ was worked out in \cite{ba0}: the fusion semiring
is identical to the one of the compact group $SO(3)$ if $n\geq 4$.
A very rough explanation for this result is the fact that $SO(3)$ is the 
(quantum) automorphism group of the $\c^*$-algebra $M_2(\C)$, and hence $SO(3)$ and $\Q_n$ have the same kind of universal property.   

Since then work on quantum permutation groups has been done  mainly in two directions.
The first one is the construction of non-classical quantum
permutation groups, i.e. quantum subgroups of $\Q_n$, using finite graphs \cite{bi1,bi2,ba2,bb,bb2,bbch}. The other direction is the study of the structure of the $\c^*$-algebra $\c(\Q_n)$. 
A general matrix representation
 was constructed in \cite{bm}, which, using the explicit description of the Haar measure
  \cite{bc,bc2}, was shown to be faithful at
 $n=4$, giving an embedding $\c(\Q_4) \subset M_4(\c(SU(2)))$ \cite{bc2}.
The case $n\geq 5$ is certainly much harder to understand, because of the non-amenability
of the discrete quantum group dual to $\Q_n$ \cite{ba0}.

A natural problem in the area is classification problem
for quantum permutation groups on $n$ points, at least for small $n$. 
This paper deals with the case $n=4$. 

The first result is a new presentation for the algebra 
$\c(\Q_4)$. This shows that $\Q_4$ is in fact isomorphic with $SO_{-1}(3)$, a $q$-analogue
of $SO(3)$ at $q=-1$. Here a few comments are in order:
$q$-deformations of classical groups have been defined by many authors at different levels of
generality (see e.g. the book  \cite{ks}) but we have never seen a non-trivial
compact quantum group $SO(3)$ at $q=-1$. For example the quantum $SO_{-1}(3)$ in \cite{po}
is just the classical $SO(3)$. 
However it would be difficult to claim
that the quantum $SO_{-1}(3)$ found here is new: the algebra $\c(SO_{-1}(3))$ is just the quotient
of $\c(SU_{-1}(3))$ by the relations making the fundamental matrix orthogonal. 

The very first consequence of the presentation result is that all
the irreducible $*$-representations of the $\c^*$-algebra $\c(\Q_4)$ are
finite dimensional, and have dimension 1, 2 or 4.
Then we show that $\c(\Q_4)$ is a deformation of $\c(SO(3))$ by an appropriate
2-cocycle. This is used in an essential way for the classification result, which is as follows.

\begin{theorem}
 The compact quantum subgroups of 
 $\Q_4$ are exactly, up to isomorphism, the following ones.
 \begin{enumerate}
  \item $\Q_4 \simeq SO_{-1}(3)$.
  \item The quantum orthogonal group $O_{-1}(2)$.
  \item $\widehat{D}_{\infty}$, the quantum dual of the infinite dihedral group.
  \item The symmetric group $S_4$ and its subgroups.
  \item The quantum group $S_4^{\tau}$, the unique non-trivial twist
  of $S_4$.
  \item The quantum group $A_5^{\tau}$, the unique non-trivial twist
  of the alternating group $A_5$.
  \item The quantum group $D_n^\tau$, $n$ even and $n\geq 6$, the unique non-trivial 
  twist of the dihedral group of order $2n$.
  \item The quantum group $DC_n^\tau$ of order $4n$, $n\geq 2$, a pseudo-twist 
  of the dicyclic group of order $4n$.
  \item The quantum group $\widehat{D}_n$, $n\geq 3$, the quantum dual 
  of the dihedral group of order $2n$.
 \end{enumerate}
 \end{theorem} 

All the quantum groups appearing in the classification are already known.
The quantum group $O_{-1}(2)$ was constructed in \cite{bi2} as the quantum automorphism group of the graph formed by two segments. The twisted quantum groups
$S_4^\tau$ and $A_5^\tau$ were constructed by Nikshych in \cite{ni}.
The twisted quantum groups $D_n^\tau$ ($n$ even and $n\geq 6$)
and pseudo-twisted $DC_n^\tau$ ($n\geq 2$) 
appear in papers by Nikshych \cite{ni} and Vainerman \cite{va}, or 
independently in Suzuki's paper \cite{su} or Masuoka's paper \cite{ma}. Note that $DC_n^\tau$ is also a pseudo-twist
of $D_{2n}$, because they have the same fusion semiring, but we prefer to refer to the dicyclic group because of the constructions in \cite{ni,ma}. At $n=2$ 
the quantum group $DC_2^\tau$ is the one corresponding to the 
$8$-dimensional historical Kac-Paljutkin example of a non-trivial Hopf algebra.
On the other hand, the occurence of  $A_5^\tau$, $D_n^\tau$ and $DC_n^\tau$
as quantum permutation groups of $4$ points seems to be new. 

The classification for finite quantum groups also uses in an essential way fundamental results of Etingof and Gelaki \cite{eg1,eg2,eg3}.

The quantum groups in the above list can be arranged in an ADE table, by using the McKay correspondence and various techniques from \cite{bbi}.

The paper is organized as follows: 2 is a preliminary section, in 3-7 we prove
the main result, and in 8-11 we write down ADE classification tables. The final section, 12, contains a few concluding remarks. 
 
\section{Compact quantum groups}

We first recall some basic facts concerning compact quantum groups
and quantum permutation groups. The book \cite{ks}
is a convenient reference for the topic of compact quantum groups,
and all the possibly missing definitions might be found there.
All the algebras are unital, and 
$\otimes$ denotes the minimal tensor product of $\c^*$-algebras
as well as the algebraic tensor product: this should cause no confusion.

\begin{definition}
A \textbf{Woronowicz algebra} is a $\c^*$-algebra $A$ endowed with  
a $*$-morphism 
$\Delta : A \to A \otimes A$ 
satisfying the coassociativity condition and the cancellation law
$$\overline{\Delta(A)(A \otimes 1)} = A \otimes A 
= \overline{\Delta(A)(1\otimes A)}$$
The morphism $\Delta$ is called
the comultiplication of $A$.
\end{definition}

The category of 
Woronowicz algebras is defined in the obvious way (see \cite{wa0} for details). 
A commutative Woronowicz algebra is necessarily isomorphic with $\c(G)$, the
function algebra on a unique compact group $G$,  
and the category of \textbf{compact quantum groups} is defined to be the category
dual to the category of Woronowicz algebras.
Hence to any Woronowicz algebra $A$ corresponds a unique compact quantum group
according to the heuristic formula $A = \c(G)$. Also to $A$
 corresponds a unique discrete quantum group $\Gamma$, with 
$A = \c^*(\Gamma)$. The quantum groups $G$ and $\Gamma$ are dual of each other.
Summarizing, we have
$$A= \c(G) = \c^*(\Gamma)$$
$$\Gamma = \widehat{G}, \quad G = \widehat{\Gamma}$$

Woronowicz's original definition for matrix compact quantum groups
\cite{wo} is still the most useful in concrete situation, and 
we have the following fundamental result \cite{wo2}.

\begin{theorem}
 Let $A$ be a $\c^*$-algebra endowed with a $*$-morphism
 $\Delta : A \to A \otimes A$. Then $A$ is a Woronowicz algebra if and only if
there exists a family of unitary matrices
$(u^\lambda)_{\lambda \in \Lambda} \in M_{d_\lambda}(A)$ satisfying the following three 
conditions.

\begin{enumerate}
\item The subalgebra $A_0$ generated by the entries $(u_{ij}^\lambda)$
of the matrices $(u^\lambda)_{\lambda \in \Lambda}$ is dense in $A$.

\item For $\lambda \in \Lambda$ and $i,j \in \{1, \ldots, d_\lambda \}$,
one has $\Delta(u_{ij}^\lambda) = \sum_{k=1}^{d_{\lambda}}
u_{ik}^\lambda \otimes u_{kj}^\lambda$.

\item For $\lambda \in \Lambda$, the transpose matrix $(u^\lambda)^t$ is
invertible.
\end{enumerate}
\end{theorem}

In fact the $*$-algebra $A_0$ in the theorem is canonically defined, and is what is now called a CQG algebra: a Hopf $*$-algebra
having all its finite-dimensional comodules equivalent to unitary ones (see \cite{ks} for details).
The counit and antipode of $A_0$, denoted respectively $\varepsilon$ and $S$, 
are refered as the counit and antipode of $A$.
The Hopf algebra $A_0$ is called the \textbf{algebra of representative functions}
on the compact quantum group $G$ dual to $A$, with another heuristic formula
$$A_0 = \r(G) = \C[\Gamma]$$

Conversely, starting from
a CQG-algebra, the universal
$\c^*$-completion yields a Woronowicz algebra in the above sense: see the book 
\cite{ks}. In fact there are possibly several different $\c^*$-norms
on $A_0$, but we will not be concerned with that problem.

As usual, a (compact) quantum subgroup $H \subset G$ corresponds to a Woronowicz algebra
surjective morphism $\c(G)\to\c(H)$, or to a surjective
Hopf $*$-algebra morphism $\r(G)\to\r(H)$.

We have the following key examples, due to Wang (respectively in \cite{wa0} and \cite{wa}).
First we need some terminology. A matrix $u \in M_n(A)$
is said to be orthogonal if $u = \overline{u}$ and 
$uu^t = I_n = u^tu$. A matrix $u$ is said to be a magic unitary if all its entries are projections, all distinct elements of a same row 
or same column are orthogonal, and sums of rows and columns are equal to 1.
A magic unitary matrix is orthogonal.

\begin{definition}
The $\c^*$-algebras $A_o(n)$ and $A_s(n)$ are constructed as follows.
\begin{enumerate}
\item $A_o(n)$ is the universal $\c^*$-algebra generated by variables $(u_{ij})_{1\leq i,j \leq n}$, with relations making $u = (u_{ij})$ an orthogonal matrix. 

\item $A_s(n)$ is the universal $\c^*$-algebra generated by variables $(u_{ij})_{1\leq i,j \leq n}$, with relations making $u = (u_{ij})$ a magic unitary matrix.
\end{enumerate} 
\end{definition}

The $\c^*$-algebras $A_o(n)$ and $A_s(n)$ are Woronowicz algebras, 
with comultiplication, counit and antipode defined by 
$$\Delta(u_{ij}) = \sum_k u_{ik}\otimes u_{kj}, \ 
\varepsilon(u_{ij}) = \delta_{ij}, \ S(u_{ij}) = u_{ji}^*$$

We now come to quantum group actions, studied e.g. in \cite{po}.
They correspond to Woronowicz algebra coactions. 

\begin{definition}
 Let $B$ be a $\c^*$-algebra. A (right) coaction of a Woronowicz
 algebra $A$ on $B$ is a $*$-homomorphism $ \alpha : B \to B \otimes A$
 satisfying the coassociativity condition and
 $$\overline{\alpha(B)(1\otimes A)} = B \otimes A$$
\end{definition}

Wang has studied quantum groups actions  on finite-dimensional $\c^*$-algebras
in \cite{wa}, where the following result is proved.

\begin{theorem}
The algebra $A_s(n)$ has the following properties.
\begin{enumerate}
\item It is the universal Woronowicz algebra coacting on $\C^n$. 
\item It is infinite-dimensional if $n \geq 4$.
\end{enumerate}  
\end{theorem}
 
 The coaction is constructed in the following manner. Let $e_1, \ldots , e_n$ be the canonical basis of $\C^n$. Then the coaction $\alpha : \C^n \to \C^n \otimes A_s(n)$
is defined by the formula
$$\alpha(e_i) = \sum_j e_j \otimes u_{ji}$$  
We refer the reader to \cite{wa} for the precise meaning of universality
in the theorem, but roughly speaking this means that the quantum group corresponding 
to $A_s(n)$, denoted $\mathcal Q_n$ (and hence $A_s(n) = \c(\Q_n))$, is the biggest one acting on $n$ points,
and deserves to be called the \textbf{quantum permutation group on $n$ points}.

Equivalently, Wang's theorem states that any Woronowicz algebra coacting faithfully on $\C^n$ is a quotient
of the Woronowicz algebra $A_s(n)$, and shows that quantum groups acting on $n$ points
(\textbf{quantum permutation groups}) 
correspond to Woronowicz algebra quotients of $A_s(n)$. 
In particular there is a Woronowicz algebra surjective morphism 
$A_s(n) \to \c(S_n)$, yielding a quantum group embedding
$S_n \subset \Q_n$. More directly, the existence
of the surjective morphism $A_s(n) \to \c(S_n)$ follows from the fact that $\c(S_n)$ is the universal
commutative $\c^*$-algebra generated by the entries of a magic unitary matrix. 
See \cite{wa} for details.

At $n=4$, by the (co)amenability result of \cite{ba0},
the study of the (compact) quantum subgroups of $\Q_4$ reduces  
to the study of 
the CQG algebra quotients of $A_s(4)_0 = \r(\Q_4)$.

\section{The quantum symmetry group of 4 points}

In this section we give  a new presentation
for the $\c^*$-algebra $A_s(4)$. The basic idea is to use an appropriate
basis of the algebra $\C^4$. In general the algebra $\C^n$ is seen as the function algebra 
on the cyclic group $\z_n$, and the Fourier transform yields
an appropriate basis for the study of many quantum permutation groups
\cite{ba2}. Here the idea is to see $\C^4$ as the function algebra on the Klein group
$\z_2\times \z_2$, and to use the Fourier transform of this group.

We will show that $A_s(4)$ is isomorphic with the following 
Woronowicz algebra.

\begin{definition} The $\c^*$-algebra
 $\c(SO_{-1}(3))$ is the $\c^*$-algebra presented by generators
$(a_{ij})_{1\leq i,j \leq 3}$ and submitted to the following relations.
\begin{enumerate}
 \item The matrix $a = (a_{ij})$ is orthogonal, 
 \item $a_{ij}a_{ik} = - a_{ik} a_{ij}$ and  $a_{ji}a_{jk} = -a_{jk}a_{ji}$, for $j \not = k$,
 \item $a_{ij} a_{kl} = a_{kl} a_{ij}$, for $ i \not = k$ and $j \not = l$,
 \item $\sum_{\sigma \in S_3} a_{1\sigma (1)}a_{2\sigma( 2)} a_{3\sigma(3)} = 1$.
\end{enumerate}
\end{definition}

It is immediate to check that $\c(SO_{-1}(3))$ is a Woronowicz algebra with
$$\Delta (a_{ij}) = \sum_k a_{ik} \otimes a_{kj} \ , \ 
\varepsilon(a_{ij}) = \delta_{ij} \ , \ S(a_{ij}) = a_{ji}$$

The last three families of relations show that, using Rosso's presentation
for the quantum group $SU_{-1}(N)$ \cite{ro}, there is a surjective
Woronowicz $\c^*$-algebra morphism $$\c(SU_{-1}(3)) \to \c(SO_{-1}(3))$$ 
The following notation will be convenient. For $i,j \in \{1,2,3\}$
with $i \not = j$, let  $\langle i, j \rangle$ be the unique  element in $\{1,2,3\}$ 
such that $\{i,j, \langle i, j \rangle\}$ = $\{1,2,3\}$.

\begin{lemma}
 Let $i,j,k,l \in \{1,2,3\}$ with $i\not=j$ and $k \not = l$. We have
 $$a_{\langle i,j \rangle \langle k,l\rangle}=
 a_{ik}a_{jl} + a_{jk}a_{il}$$  
\end{lemma}
 
\begin{proof}
For $i,j \in \{1,2,3\}$, let $i_1,i_2,j_1,j_2$ be the elements of $\{1,2,3\}$ satisfying
$\{i,i_1,i_2\} = \{1,2,3\} = \{j,j_1,j_2\}$.
The well-known formula for the antipode in $\c(SU_{-1}(3))$ gives
$$S(a_{ji})= a_{i_{1}j_{1}}a_{i_{2}j_{2}} + a_{i_{2}j_{1}} a_{i_{1}j_{2}}$$
But since the matrix $a$ is orthogonal, we have $a_{ij}=S(a_{ij})$, which proves the statement.  
\end{proof}

We now prove the presentation result.

\begin{theorem}
 We have a Woronowicz algebra isomorphism $\c(\Q_4) \simeq \c(SO_{-1}(3))$.
\end{theorem}
 
\begin{proof}
Let us show that  $\c(SO_{-1}(3))$ is the universal Woronowicz algebra 
coacting on $\C^4$. Let $e_1, e_2, e_3, e_4$ be the canonical basis of $\C^4$. We use the following basis
\begin{align*}
1 & = e_1+e_2+e_3+e_4\\
\varepsilon_1 & = e_1-e_2-e_3+e_4\\
\varepsilon_2 & = e_1-e_2+e_3-e_4\\
\varepsilon_3 &= e_1+e_2-e_3-e_4 = \varepsilon_1\varepsilon_2
\end{align*}
This basis is just obtained by using 
the Fourier transform of the group $\z_2 \times \z_2$.
We have $\varepsilon_i^2 = 1$, $\varepsilon_i \varepsilon_j = \varepsilon_j \varepsilon_i
= \varepsilon_{\langle i, j \rangle}$ for $i \neq j$,
$\varepsilon_i^*= \varepsilon_i$, and these relations define a presentation of $\C^4$. 
Consider the linear map $\alpha : \C^4 \to \C^4 \otimes \c(SO_{-1}(3))$ defined
by $\alpha(1) = 1 \otimes 1$ and $\alpha(\varepsilon_i) = \sum_j \varepsilon_j \otimes a_{ji}$.
It is clear that $\alpha$ is coassociative, and it remains to check that $\alpha$ is 
a $*$-algebra morphism. We have
\begin{align*}
 \alpha(\varepsilon_i)^2 & = \sum_{k,l} \varepsilon_k \varepsilon_l
 \otimes a_{ki}a_{li}\\
 & = \sum_{k} \varepsilon_k^2 \otimes a_{ki}^2 + \sum_{k \neq l} \varepsilon_k \varepsilon_l \otimes
 a_{ki}a_{li}\\
 & = 1 \otimes (\sum_{k} a_{ki}^2) + \sum_{k<l} \varepsilon_k\varepsilon_l \otimes (a_{ki}a_{li} + 
 a_{li}a_{ki}) \\
 & = 1 \otimes 1 = \alpha(\varepsilon_i^2)
\end{align*}
By using Lemma 3.1, we have, for $i \not = j$,
\begin{align*}
 \alpha(\varepsilon_i)\alpha(\varepsilon_j) & = \sum_{k,l} \varepsilon_k \varepsilon_l
 \otimes a_{ki}a_{lj}\\
 & = \sum_{k} \varepsilon_k^2 \otimes a_{ki}a_{kj} + \sum_{k \neq l} \varepsilon_k \varepsilon_l \otimes
 a_{ki}a_{lj}\\
 & = 1 \otimes (\sum_{k} a_{ki}a_{kj}) + \sum_{k<l} \varepsilon_k\varepsilon_l \otimes (a_{ki}a_{lj} + 
 a_{li}a_{kj})\\
 & = \sum_{k<l} \varepsilon_{\langle k,l\rangle} \otimes a_{\langle k, l\rangle \langle i, j\rangle}\\
 & = \sum_k \varepsilon_k \otimes a_{k \langle i, j\rangle} =
  \alpha(\varepsilon_{\langle i,j\rangle})
\end{align*}
and since $\alpha(\varepsilon_i^*) = \alpha(\varepsilon_i)$, we conclude that $\alpha$ is a coaction.

Consider now a Woronowicz algebra $A$ coacting on $\C^4$, with coaction 
$\beta : \C^4 \to \C^4 \otimes A$. We put $\varepsilon_0 =1$ and we can write
$$\beta(\varepsilon_i) = \sum_j \varepsilon_j \otimes x_{ji}$$
with $x_{j0} = \delta_{j0}$. We have $x_{ij}^* =x_{ij}$ since $\beta$ is
a $*$-morphism. 
Let $\phi : \C^4 \to \C$ be the classical normalized integration
map: $\phi(e_i) =1/4$, and hence $\phi(\varepsilon_i) = 0$ if $i>0$. 
The linear map $\phi$ is $A$-colinear, and thus the linear map 
$\C^4 \otimes \C^4 \to \C$, $x \otimes y \mapsto \phi(xy)$, is $A$-colinear.
It follows that the matrix $x = (x_{ij})$ is orthogonal, and hence
$$x_{01}^2 + x_{02}^2 + x_{03}^2 = 0 =  x_{01}x_{01}^* + x_{02}x_{02}^* + x_{03}x_{03}^*$$
Therefore $x_{0i} = \delta_{0i}$, the matrix $x' =(x_{ij})_{1\leq i,j\leq3}$ is orthogonal and
 $$\beta(\varepsilon_i) = \sum_{j=1}^3 \varepsilon_j \otimes x_{ji}, \ 1 \leq i \leq 3.$$
We have 
\begin{align*}
 1\otimes 1 & =  \beta(1) = \beta(\varepsilon_i^2) = \beta(\varepsilon_i)^2 \\
  & = \sum_{k,l} \varepsilon_k \varepsilon_l
 \otimes x_{ki}x_{li}\\
 & = \sum_{k} \varepsilon_k^2 \otimes x_{ki}^2 + \sum_{k \neq l} \varepsilon_k \varepsilon_l \otimes
 x_{ki}x_{li}\\
 & = 1 \otimes (\sum_{k} x_{ki}^2) + \sum_{k<l} \varepsilon_k\varepsilon_l \otimes (x_{ki}x_{li} + 
 x_{li}x_{ki}) \\
 & = 1 \otimes 1 + \sum_{k<l} \varepsilon_k\varepsilon_l \otimes (x_{ki}x_{li} + 
 x_{li}x_{ki})
\end{align*}
 and hence for $k \not = l$, we have $x_{ki}x_{li} = - x_{li}x_{ki}$, and using the antipode
 ($S(x_{ij} ) = x_{ji}$ since the matrix $x'$ is orthogonal), we also get 
 $x_{ik}x_{il} = - x_{il}x_{ik}$.
 For $i \not = j$, we have
 \begin{align*}
 \alpha(\varepsilon_{\langle i, j\rangle}) &= 
 \alpha(\varepsilon_i \varepsilon_j)  = \alpha(\varepsilon_i)\alpha(\varepsilon_j) = \sum_{k,l} \varepsilon_k \varepsilon_l
 \otimes x_{ki}x_{lj}\\
 & = \sum_{k} \varepsilon_k^2 \otimes x_{ki}x_{kj} + \sum_{k \neq l} \varepsilon_k \varepsilon_l \otimes
 x_{ki}x_{lj}\\
 & = 1 \otimes (\sum_{k} x_{ki}x_{kj}) + \sum_{k<l} \varepsilon_k\varepsilon_l \otimes (x_{ki}x_{lj} + 
 x_{li}x_{kj})\\
 & = \sum_{k<l} \varepsilon_{\langle k,l\rangle} \otimes (x_{ki}x_{lj}+ x_{li}x_{kj})
\end{align*}
 and hence $x_{\langle k,l \rangle \langle i , j \rangle} = x_{ki}x_{lj}+ x_{li}x_{kj}$. Similarly,
 since $\varepsilon_i \varepsilon_j = \varepsilon_j \varepsilon_i$, we have
  $x_{\langle k,l \rangle \langle i , j \rangle} = x_{kj}x_{li}+ x_{lj}x_{ki}$.
  By combining these two relations, we get, for $i \not =j$ and $k \not = l$,
  $$[x_{ki},x_{lj}] = [x_{kj},x_{li}]$$
  Now using the antipode, 
  we also have
  $$[x_{ki},x_{lj}] = [x_{li},x_{kj}]$$
  This finally leads to 
  $$[x_{ki},x_{lj}] = [x_{li},x_{kj}] =[x_{kj},x_{li}] = -[x_{li},x_{kj}]$$
  and we conclude that $x_{li}x_{kj} = x_{kj}x_{li}$.  
  The quantum determinant relation now follows.
 \begin{align*}
  & x_{11}x_{22}x_{33} + x_{11}x_{23}x_{32} + x_{12}x_{21}x_{33} + x_{12}x_{23}x_{31} +
  x_{13}x_{22}x_{31} + x_{13}x_{23}x_{31} \\
  = & x_{11}(x_{22}x_{33} +x_{23}x_{32}) + x_{12}(x_{21}x_{33} + x_{23}x_{31}) +
  x_{13}(x_{22}x_{31} + x_{23}x_{31}) \\
  = & x_{11}^2 + x_{12}^2 + x_{13}^2 =1
 \end{align*}
Therefore we get a Woronowicz algebra morphism 
$\c(SO_{-1}(3)) \to A$ commuting with the respective coactions,
and we are done.
\end{proof}
 
\noindent
\textbf{Remark}. The proof of the theorem furnishes the following
concrete isomorphism $\c(SO_{-1}(3)) \simeq \c(\Q_4)$:
$$
\left(\begin{array}{cccc}
1 & 0 & 0 & 0 \\ 
0 & a_{11} & a_{12} & a_{13} \\ 
0 & a_{21} & a_{22} & a_{23} \\ 
0 & a_{31} & a_{32} & a_{33}
\end{array}\right)
\longmapsto
\frac{1}{4}
M
\left(
\begin{array}{llll}
u_{11} & u_{12} & u_{13} & u_{14} \\ 
u_{21} & u_{22} & u_{23} & u_{24} \\ 
u_{31} & u_{32} & u_{33} & u_{34} \\ 
u_{41} & u_{42} & u_{43} & u_{44}
\end{array}\right) M
$$ 
where 
$$M=\left(
\begin{array}{cccc}
1 & 1 & 1 & 1 \\ 
1 & -1 & -1 & 1 \\ 
1 & -1 & 1 & -1 \\ 
1 & 1 & -1 & -1
\end{array}\right)$$
 The very first consequence of the isomorphism theorem is
 the following result.
 
 \begin{proposition}
  Any irreducible Hilbert space representation of the $\c^*$-algebra
  $A_s(4)=\c(\Q_4)$ is finite dimensional, and has dimension $1$, $2$ or $4$.
 \end{proposition}
 
 \begin{proof}
 We just have to prove the result for $\c(SO_{-1}(3))$. Let 
 $$\pi :   \c(SO_{-1}(3)) \to \mathcal B(H)$$
 be an irreducible representation on a Hilbert space $H$. For simplicity we just write
 $A = \c(SO_{-1}(3))$. The elements $a_{ij}^2$ belong to the center of $A$, hence
 $\pi(a_{ij}^2)$ is a scalar. Using this fact and the commutation and anticommutation
 relations in $A$, we see that $\pi(A)$ is finite dimensional with 
 $\dim(\pi(A)) \leq 2^9$. Il follows that $H$ is finite-dimensional. The relations
 of Lemma 3.1
 \begin{align*}
  a_{13} & = a_{21}a_{32} + a_{22}a_{31} \\
  a_{23} & = a_{11}a_{32} + a_{12}a_{31} \\
  a_{33} & = a_{11}a_{22} + a_{12}a_{21}
 \end{align*}
now ensure that $\dim(\pi(A)) \leq 2^6$.
On each column of the orthogonal matrix $(\pi(a_{ij}))$, there is at least 
a nonzero element, which is invertible by the irreducibility of $\pi$, and 
from this we see that $\dim(\pi(A)) \leq 2^5$. If $n= \dim(H)$, we have 
$n^2 = \dim(\pi(A))$, and hence $n \leq \sqrt{2}\cdot 4 < 6$.
We also see from the anticommutation relations that if $n \not = 1$, then $n$
is necessarily even: this completes the proof.  
 \end{proof}

\section{Twistings}

We discuss here the 2-cocycle deformation 
procedure, that we call twisting, at a purely algebraic level.
This method was initiated by Drinfeld, and studied in a systematic way 
in the dual framework by Doi \cite{do}. 

Let $H$ be a Hopf algebra. We use Sweedler's notation
$\Delta(x) = x_{1} \otimes x_{2}$. Recall (see e.g. \cite{do})
that a 2-cocycle is a convolution invertible linear map
$\sigma : H \otimes H \longrightarrow \c$ satisfying
$$\sigma(x_{1}, x_{2}) \sigma(x_{2}y_{2},z) =
\sigma(y_{1},z_{1}) \sigma(x,y_{2} z_{2})$$
and $\sigma(x,1) = \sigma(1,x) = \varepsilon(x)$, for $x,y,z \in H$.

Following \cite{do} and \cite{sc1}, we associate various algebras to
a 2-cocycle. First consider the algebra 
$_{\sigma} \! H$. As a vector space we have $_{\sigma} \! H = H$ and the product
of $_{\sigma}H$ is defined to be
$$\{x\}  \{y\} = \sigma(x_{1}, y_{1}) \{x_{2} y_{2}\}, 
\quad x,y \in H,$$
where an element $x \in H$ is denoted $\{x\}$, when viewed as an element 
of $_{\sigma} \!H$.

We also have the algebra $H_{\sigma^{-1}}$, where $\sigma^{-1}$ denotes the convolution inverse of $\sigma$. As a vector space we have
$H_{\sigma^{-1}} = H$ and the product of 
$H_{\sigma^{-1}}$ is defined to be
$$\langle x \rangle \langle y \rangle = \sigma^{-1}(x_{2}, y_{2}) \langle x_{1} y_{1} \rangle, 
\quad x,y \in H.$$
where an element $x \in H$ is denoted $\langle x \rangle$, when viewed as an element 
of $H_{\sigma^{-1}}$. The cocycle condition ensures that $_{\sigma} \! H$
and $H_{\sigma^{-1}}$ are associative algebras with $1$ as a unit.

Finally we have the 
Hopf algebra $H^{\sigma} = {_{\sigma} \! H}\!_{\sigma^{-1}}$.  
As a coalgebra $H^{\sigma} 
= H$. The product of $H^{\sigma}$ is defined to be
$$[x] [y]= \sigma(x_{1}, y_{1})
\sigma^{-1}(x_{3}, y_{3}) [x_{2} y_{2}], 
\quad x,y \in H,$$
where an element $x \in H$ is denoted $[x]$, when viewed as an element 
of $H^{\sigma}$, 
and we have the following formula for the antipode of 
$H^{\sigma}$:
$$S^{\sigma}([x]) = \sigma(x_{1}, S(x_{2}))
\sigma^{-1}(S(x_{4}), x_{5})[ S(x_{3})].$$
The Hopf algebras $H$ and $H^{\sigma}$ have equivalent tensor categories of comodules
\cite{sc1}.

Very often 2-cocyles are induced by simpler quotient Hopf algebras (quantum subgroups).
More precisely let $\pi : H \to K$ be a Hopf algebra surjection and let
$\sigma : K \otimes K \to \C$ be a 2-cocycle on $K$.
Then $\sigma_{\pi} = \sigma \circ (\pi \otimes \pi) : H \otimes H \to \C$ is a 2-cocycle.

\begin{lemma}
Let $\pi : H \to K$ be a Hopf algebra surjection and let
$\sigma : K \otimes K \to \C$ be a $2$-cocycle on $K$. We have an injective
algebra map
\begin{align*}
\theta : H^{\sigma_\pi} & \longrightarrow {_\sigma \! K} \otimes H \otimes K_{\sigma^{-1}} \\
[x] & \longmapsto \{\pi(x_1)\} \otimes x_{2} \otimes \langle \pi(x_3) \rangle 
\end{align*}
\end{lemma}

\begin{proof}
It is a direct verification that $\theta$ is an algebra map, and the injectivity follows by using the counit of $K$, available as a linear map. 
\end{proof}

This lemma, when $K$ is finite dimensional and $H$ is commutative,
furnishes a faithful representation of $H^{\sigma_{\pi}}$.

The following lemma will be used in the section concerning the classification.

\begin{lemma}
 Let $\pi_1,\pi_2 : H \to K$ be two Hopf algebra surjections and let
$\sigma : K \otimes K \to \C$ be a $2$-cocycle on $K$.
Assume that there exists  a Hopf algebra automorphism $u : H \to H$
such that $\pi_2 \circ u = \pi_1$. Then the Hopf algebra $H^{\sigma_{\pi_1}}$ 
and $H^{\sigma_{\pi_2}}$ are isomorphic.
\end{lemma}

\begin{proof}
 The announced isomorphism is the map $H^{\sigma_{\pi_1}} \to H^{\sigma_{\pi_2}}$, $[x] \mapsto [u(x)]$. 
\end{proof}



The following lemma is useful to classify certain quotient Hopf algebras.

\begin{lemma}
Let $\pi : H \to K$ be a Hopf algebra surjection and let 
 $\sigma : K \otimes K \to \C$ be a $2$-cocycle.
 Then there is a bijection between the following data.
\begin{enumerate}
 \item Surjective Hopf algebra maps $ f : H \to L$ such that there exists
 a Hopf algebra map $g : L \to K$ satisfying $g \circ f = \pi$.
 \item Surjective Hopf algebra maps $ f' : H^{\sigma_\pi} \to L'$ such that there exists
 a Hopf algebra map $g' : L' \to K^{\sigma}$ satisfying $g' \circ f' = [\pi]$.
\end{enumerate}
Assume moreover that $H$ and $K$ are Hopf $*$-algebras, that $\pi$ is a Hopf $*$-algebra map and that $H^{\sigma_\pi}$ and $K^\sigma$ admit Hopf $*$-algebra structures
induced by those of $H$ and $K$ respectively (i.e. $[x]^* = [x^*]$).
Then the above correspondence remains true for Hopf $*$-algebras and 
Hopf $*$-algebra surjections.
\end{lemma}

\begin{proof}
 Let us start with surjective Hopf algebra maps $ f : H \to L$ and 
 $g : L \to K$ satisfying $g \circ f = \pi$. Then the Hopf algebra maps 
 $ [f] : H^{\sigma_\pi} \to L^{\sigma_g}$ and $[g] : L^{\sigma_g} \to K^{\sigma}$
 satisfy $[g] \circ [f] = [\pi]$. We get the other side of the correspondence
 by using the cocycle $\sigma^{-1}$. The Hopf $*$-algebra assertion is then immediate.
\end{proof}

We conclude the section by mentioning the pseudo-twisting procedure, considered 
in the dual framework in \cite{ev,va,ni}. Consider a convolution invertible
map $\sigma : H \otimes H \to \C$. If $\sigma$ is a pseudo-2-cocycle, it is possible to define  a Hopf algebra $H^\sigma$ as above, and $H^\sigma$ has the same fusion semiring as $H$.
Some examples of pseudo-twisted Hopf algebras will occur in Section 7.

\section{Twisting and quantum $SO_{-1}(3)$}

We now proceed to show that the quantum group $SO_{-1}(3)$
of the previous section may be obtained by twisting the group $SO(3)$.   

Let $H$ be the algebra of representative functions
on the compact group $SO(3)$, with the canonical 
coordinate functions denoted  $x_{ij}$. We consider the Klein group, denoted
$\V$ as usual,
$$\V = \z_2 \times \z_2 = \langle t_1 , t_2 \ | t_1^2 = t_2^2 =1, \ t_1 t_2 = t_2 t_1\rangle$$
and we put $t_3 = t_1t_2$. The restriction of functions on $SO(3)$ to its diagonal
subgroup gives a Hopf algebra surjection
\begin{align*}
\pi_d : H & \longrightarrow \c[\V]\\
  x_{ij} & \longmapsto \delta_{ij} t_i
\end{align*}
Let $\sigma : \V \times \V \longrightarrow \C^*$ be the unique
bimultiplicative map such that $\sigma(t_i,t_j) = -1$ if $i \leq j$
and $\sigma(t_i,t_j) = 1$ otherwise.
Then $\sigma$ is a usual group 2-cocycle, and its unique linear extension to the group algebra
$\C[\V]\otimes \C[\V]$ is a 2-cocycle in the previous sense, still denoted $\sigma$.
We get a cocycle $\sigma_d = \sigma \circ (\pi_d \otimes \pi_d)$ on $H$.
Note that $\sigma_d = \sigma_d^{-1}$.

\begin{theorem}
 The Hopf algebra $H^{\sigma_d}$ is a CQG-algebra, and the Woronowicz 
 algebras $\c^*(H^{\sigma_d})$ and $\c(SO_{-1}(3))$ are isomorphic.
\end{theorem}

\begin{proof} We have
 \begin{align*}
  [x_{ij}][x_{kl}] & = \sum_{p,q,r,s} \sigma_d(x_{ip},x_{kq})
  \sigma_d^{-1}(x_{rj},x_{sl})[x_{pr}x_{qs}]\\
  & = \sigma(t_i,t_k)\sigma(t_j,t_l)[x_{ij}x_{kl}] 
 \end{align*}
 and more generally
 $$[x_{i_1j_1}] [x_{i_2j_2}] \ldots [x_{i_rj_r}] = 
 \left(\prod_{p<q}\sigma(t_{i_p},t_{i_q})\right)
 \left(\prod_{p<q}\sigma(t_{j_p},t_{j_q})\right)
 [x_{i_1j_1} x_{i_2j_2} \ldots x_{i_rj_r}]$$
Thus we see
that there is a Hopf algebra morphism $\r(SO_{-1}(3)) \to H^{\sigma_d}$ mapping  
$a_{ij}$ to $[x_{ij}]$. This morphism is clearly surjective.
Now combining Theorem 3.1 and \cite{ba0},  
we see that the representation semiring of $\r(SO_{-1}(3))$ is the same as the one of $SO(3)$,
and thus by a standard argument the two cosemisimple
Hopf algebras $\r(SO_{-1}(3))$ and $H^{\sigma_d}$ are isomorphic. The theorem follows. 
\end{proof}

We now  recover the faithful representation of \cite{bc2}, under a different
form. The same kind of embedding was constructed for $SU_{-1}(2)$ in \cite{z}.
It will be used in the next section. First we need the following lemma, whose
proof is immediate.

\begin{lemma}
 Let $H$ be a cosemisimple Hopf algebra with Haar measure $h : H \to \C$, and let 
 $\pi : H \to K$ be a Hopf algebra morphism. Let 
 $\varphi : K \to \C$ be a linear map such that $\varphi(1)=1$. Then for $x \in H$, we have
 $$h(x) =\varphi(\pi(x_1)) \varphi(\pi(x_3))h(x_2)$$
\end{lemma}

We use the Pauli matrices
$$\tau_1 = 
\left(\begin{array}{ll}
i & 0 \\ 
0 & -i 
  \end{array}\right),
\quad
\tau_2 =\left(\begin{array}{ll}
0 & 1 \\ 
-1 & 0 
  \end{array}\right),
\quad \tau_3 = \tau_2\tau_1 =
\left(\begin{array}{ll}
0 & -i \\ 
-i & 0 
  \end{array}\right)$$

\begin{theorem}
 We have a $\c^*$-algebra embedding
 \begin{align*}
 \theta : \c(SO_{-1}(3)) & \longrightarrow M_2(\C) \otimes M_2(\C) \otimes \c(SO(3)) \\ 
 a_{ij} & \longmapsto \tau_i \otimes \tau_j \otimes x_{ij}
 \end{align*}
\end{theorem}

\begin{proof}
 The first thing to note is that the twisted group
 algebra $\C[\V]_\sigma = \C_\sigma[\V]$ is isomorphic to
 $M_2(\C)$ via $\{t_i\} \mapsto \tau_i$.
 Then combining this with Theorem 5.1 and Lemma 4.1, we get a $*$-algebra
 embedding $f : \r(SO_{-1 }(3)) \to M_2(\C) \otimes M_2(\C) \otimes H$, inducing
 the announced $*$-algebra map. 
 Let ${\rm tr}$ be the normalized trace on $M_2(\C)$, and let $h$ be the Haar mesure on $H$.
 Then $({\rm tr} \otimes {\rm tr} \otimes h)\circ f$ is,  
 by the invariance of the Haar measure by cocycle twisting and Lemma 5.1, the Haar measure on $\r(SO_{-1}(3))$. We can conclude as in \cite{bc2}, using the amenability
 of $\widehat{\Q_4}$, that the $*$-algebra morphism
 $\c(SO_{-1}(3))  \longrightarrow M_2(\C) \otimes M_2(\C) \otimes \c(SO(3))$
 is injective.
\end{proof}

\section{Quantum subgroups of $\Q_4$}

We now study the quantum subgroups of $\Q_4$. They correspond
to Woronowicz algebra quotients of $\c(\Q_4)$ and hence of $\c(SO_{-1}(3))$, 
or equivalently to Hopf
$*$-algebra quotients of $\r(SO_{-1}(3))$. We work at the Hopf algebra level. The strategy for proving Theorem 1.1 is the following one.
Our first goal is the classification of the quantum subgroups of 
$SO_{-1}(3)$ diagonally containing the Klein subgroup (see the notation below).
This is achieved in Theorem 6.1, by using the technical results
6.1-6.5. Then another technical result, Lemma 6.6, shows that the general classification
follows from the classification follows from the classification 
for some proper quantum subgroups. The classification is then performed for these 
quantum subgroups, with the exception of the quantum group $O_{-1}(2)$, that
will be treated in the next section.

First we need some notations. Most of them have been introduced in the previous section.
\begin{enumerate}
\item The Hopf $*$-algebra $H$ is the Hopf algebra $\r(SO(3))$, with the coordinate functions
denoted $x_{ij}$. We denote by $I_d$ the ideal of $H$ generated by the elements $x_{ij}$, $i \not = j$.

\item The Klein group is the group $\V = \z_2 \times \z_2 = \langle t_1 , t_2 \ | t_1^2 = t_2^2 =1, \ t_1 t_2 = t_2 t_1\rangle$, and we put $t_3 = t_1t_2$. 

\item The group of diagonal matrices in $SO(3)$ is denoted $D$. It is isomorphic
with $\V$ via the injective group morphism
$$t_1 \longmapsto \left(
\begin{array}{lll}
1 & 0 & 0 \\ 
0 & -1 & 0 \\ 
0  & 0 & -1
\end{array}\right), \quad
t_2 \longmapsto 
\left(
\begin{array}{lll}
-1 & 0 & 0 \\ 
0 & 1 & 0 \\ 
0  & 0 & -1
\end{array}\right)$$
   
\item We have a ``diagonal'' Hopf $*$-algebra surjection $\pi_d : H \to \C[\V]$, $x_{ij} \mapsto
\delta_{ij} t_i$, induced by the previous diagonal embedding $\V \hookrightarrow SO(3)$.
It is clear that $I_d = \ker(\pi_d)$.

\item The 2-cocycle on $\V$ defined in the previous section is denoted $\sigma$, and it induces a 
2-cocycle $\sigma_d$ on $H$. More generally, if $K$ is a Hopf algebra and $f : K \to \C[\V]$ is a 
Hopf algebra map, then $\sigma_f = \sigma \circ (f \otimes f)$ is a 2-cocycle on $K$.

\item The Hopf $*$-algebra $H'$ is the Hopf algebra $\r(SO_{-1}(3))$. Thanks to Theorem 5.1
in the previous section, we identify it with $H^{\sigma_d}$.
We have a diagonal surjective Hopf $*$-algebra map $\pi_d' : H' \to \C[\V]$, $a_{ij} \mapsto
\delta_{ij}t_i$. 
$\pi_d'$ induces a quantum group embedding $\V \subset SO_{-1}(3)$, called the diagonal embedding
of $\V$ into $SO_{-1}(3)$.
In fact $\pi_d'$ is just $[\pi_d]$, after the identification $H' = H^{\sigma_d}$. We put $I_d' = \ker(\pi_d')$. Clearly this is the ideal generated by the elements $a_{ij}$, $i \not = j$.
\end{enumerate}

\subsection{Quantum subgroups diagonally containing the Klein subgroup.}

We begin with a lemma.

\begin{lemma}
Let $X \subset SO(3)$ be a compact subgroup with $D \subset X$, and let $g : \r(X) \to \C[\V]$
be the corresponding Hopf $*$-algebra surjection. Then $\r(X)^{\sigma_g}$ is a quotient Hopf $*$-algebra of
$\r(SO_{-1}(3))$.

Conversely,  
 let $G\subset SO_{-1}(3)$ be a compact quantum subgroup, and let  $f : H' \to \r(G)$ be the  corresponding surjective
 Hopf $*$-algebra map. Assume that there exists a Hopf $*$-algebra map
 $g' : \r(G) \to \C[\V]$ such that $g' \circ f = \pi_d'$: the corresponding composition of embeddings
 $\V \subset G \subset SO_{-1}(3)$ is the diagonal embedding.
 Then there exists a compact subgroup $X$ of $SO(3)$ with
 $D \subset X \subset SO(3)$, with a Hopf $*$-algebra isomorphism
 $$\r(G) \simeq \r(X)^{\omega}$$
for  a $2$-cocyle $\omega$ induced by $\sigma$. 
\end{lemma}

\begin{proof} The Hopf $*$-algebra quotients of 
$H= \r(SO(3))$ all are isomorphic with $\r(X)$, for a compact subgroup
$X \subset SO(3)$. Hence this lemma is just a reformulation of Lemma 4.3 annd its proof.
\end{proof}

Using the well-known list of compact subgroups of $SO(3)$, we will get the first classification result.
The subgroups containing $D$ are: $SO(3)$, $O(2)$, $S_4$, $A_4$, $D_4$, $D\simeq \z_2\times \z_2$,
$D_n$ with $n>4$ even and $A_5$.

We need the following lemma, which collects various results on 2-cocycles 
deformations in the litterature.

\begin{lemma}
 The groups $A_4$, $D_4$ and  $\z_2\times \z_2$ do not have any non-trivial twist,
 and the groups $S_4$, $A_5$ and $D_n$ for $n>4$ even, have exactly one 
 non-trivial twist. 
\end{lemma}

\begin{proof}
The assertion for the dihedral groups is due to Masuoka \cite{mas}, and for $\z_2 \times \z_2$
this is trivial. 
 For the other groups  we have the following general fact, due to Etingof-Gelaki \cite{eg1,eg3}
 and Davydov \cite{da}. If $G$ is a finite group, any non-trivial twist of
 $\c(G)$ arises from the following data: a non-normal subgroup $K$ endowed
 with a (necessarly non trivial) 2-cocycle $\alpha$ on $K$ such that the twisted group
 algebra of $\C_\alpha[K]$ is a matrix algebra. For the groups $A_4$, $S_4$ and $A_5$, the only 
 possible candidate for $K$ is a Klein subgroup $\V$, which is normal (and hence unique) for $A_4$, and
 hence we have the assertion for $A_4$. Also there is only one non trivial cohomology class
  in $H^2(\V, \C^*)$, and since the non-normal Klein subgroups of these groups
  are all conjugate, by Lemma 4.2 there is at most one possible twist for $S_4$ and $A_5$. These twists were constructed by Nikshych \cite{ni}.
\end{proof}

\begin{lemma}
 The compact groups $SO(3)$, $O(2)$,
$D_n$ with $n$ even and $n > 4$, and $A_5$ are not quantum subgroups of $SO_{-1}(3)$.
\end{lemma}

\begin{proof}
This follows from the fact that  $\c(S_4)$ is the maximal abelian quotient
of $\r(\Q_4) \simeq \ \r(SO_{-1}(3))$. 
\end{proof}

\begin{lemma}
The finite quantum groups $D_n^\tau$, $n$ even and $n>4$, 
 $S_4^{\tau}$,  and 
  $A_5^{\tau}$, the respective unique non-trivial twists of $D_n$, $S_4$ and $A_5$,
  all occur as quantum subgroups of $SO_{-1}(3)$ diagonally containing $\V$.
\end{lemma}

\begin{proof}
 We know from Lemma 6.1 that for a compact group $X$ with $D\subset X \subset SO(3)$, 
 there exists a
 2-cocycle $\omega$ on $\r(X)$ such that $\r(X)^{\omega}$ is a quotient Hopf $*$-algebra
 of   $\r(SO_{-1}(3))$. Using  the previous lemma, we see that $D_n^\tau$ anf $A_5^\tau$
 are quantum subgroups of $SO_{-1}(3)$ diagonally containing $\V$. The quantum group $S_4^{\tau}$ is realized 
 as a quantum subgroup of $\Q_4$ in \cite{bi0}.
 Still denoting by $u_{ij}$ the generators of $\c(S_4^\tau)$, we have, 
by using the presentation of $\c(S_4^\tau)$ given in \cite{bi0}, 
 a surjective Hopf $*$-algebra map $\c(S_4^\tau) \to \C[\V]$
 $$
 \left(
\begin{array}{llll}
u_{11} & u_{12} & u_{13} & u_{14} \\ 
u_{21} & u_{22} & u_{23} & u_{24} \\ 
u_{31} & u_{32} & u_{33} & u_{34} \\ 
u_{41} & u_{42} & u_{43} & u_{44}
\end{array}\right) \longmapsto
\left(
\begin{array}{llll}
\delta_{1}   & \delta_{t_1} & \delta_{t_2} & \delta_{t_3} \\ 
\delta_{t_1} & \delta_{1}   & \delta_{t_3} & \delta_{t_2} \\ 
\delta_{t_2} & \delta_{t_3} & \delta_{1}   & \delta_{t_1} \\ 
\delta_{t_3} & \delta_{t_2} & \delta_{t_1}  & \delta_{1}
\end{array}\right)$$
Composed with the canonical morphisms $\r(SO_{-1}(3)) \simeq\r(\Q_4) \to C(S_4^\tau)$, 
this gives embeddings $V \subset S_4^\tau \subset SO_{-1}(3)$, and the composition is the diagonal embedding.  
\end{proof}

We will also need the quantum group $O_{-1}(2)$, also considered in \cite{bbc}.

\begin{definition} The $*$-algebra
 $\r(O_{-1}(2)))$ is the $*$-algebra presented by generators
$(a_{ij})_{1\leq i,j \leq 2}$ and submitted to the following relations.
\begin{enumerate}
 \item The matrix $a = (a_{ij})$ is orthogonal, 
 \item $a_{ij}a_{ik} = - a_{ik} a_{ij}$ and  $a_{ji}a_{jk} = -a_{jk}a_{ji}$, for $j \not = k$,
 \item $a_{ij} a_{kl} = a_{kl} a_{ij}$, for $ i \not = k$ and $j \not = l$,
\end{enumerate}
\end{definition}

It is immediate to check that $\r(O_{-1}(2))$ is a Hopf $*$-algebra with
$$\Delta (a_{ij}) = \sum_k a_{ik} \otimes a_{kj} \ , \ 
\varepsilon(a_{ij}) = \delta_{ij} \ , \ S(a_{ij}) = a_{ji}$$
and is a CQG-algebra. Also it is easy to construct a surjective 
Hopf $*$-algebra map $\r(SO_{-1}(3)) \to \r(O_{-1}(2))$, yielding an embedding
$O_{-1}(2) \subset SO_{-1}(3)$ such that the composition of the embeddings 
$\V \subset O_{-1 }(2) \subset SO_{-1}(3)$ is the diagonal one. 

\begin{lemma}
 The Hopf $*$-algebra $\r(O_{-1}(2))$ is a twist of $\r(O(2))$, by a $2$-cocycle
 induced by the diagonal Klein subgroup of $O(2)$. Moreover this is the only non-trivial 
 such twist of $\r(O(2))$.
\end{lemma}

\begin{proof}
The first assertion is proved in \cite{bbc}
where it is shown that $O_{-1}(2)$ is the quantum symmetry group of the square, 
see also \cite{bi2}.
All the  Klein subgroups of $O(2)$ are conjugate, and we have the uniqueness result using
Lemma 4.2.
\end{proof}

Combining the above results,  we can state the first classification result.

\begin{theorem}
 The  compact quantum subgroups $G$ of Ŝ$SO_{-1}(3)$ with $\V \subset G$, and such that
 the composition of the embeddings $\V \subset G \subset SO_{-1}(3)$ is the 
 diagonal embedding, are exactly, up to isomorphism, the following ones.
 \begin{enumerate}
  \item $SO_{-1}(3)$.
  \item $O_{-1}(2)$.
  \item $\z_2 \times \z_2$, $D_4$, $A_4$ or $S_4$.
  \item $D_n^\tau$, $n$ even and $n>4$, the only non-trivial twist
  of the dihedral group $D_n$.
  \item $S_4^{\tau}$, the only non-trivial twist
  of the symmetric group $S_4$.
  \item $A_5^{\tau}$, the only non-trivial twist
  of the alternating group $A_5$.
 \end{enumerate}
\end{theorem}

\subsection{The general case.}

We now deal with general quantum subgroups. 
So let us consider a Hopf $*$-algebra surjection 
$f : H' = \r(SO_{-1}(3))\to L$ and let $J=\ker f$.
Consider the Hopf $*$-algebra $\tilde{L} = L/(J \cap I_d')$. 
It is clear that $\tilde{L} = \r(G)$ for one of the quantum groups of Theorem 6.1.
The following lemma ensures that if $\tilde{L}$ is $\r(SO_{-1}(3))$, then already
$L$ was $ \r(SO_{-1}(3))$.  

\begin{lemma}
 Let $J$ be a Hopf $*$-ideal in $H'$. If $J \cap I_d' = (0)$, then 
 $J=(0)$.
\end{lemma}

\begin{proof}
Consider the embedding 
\begin{align*}
 \theta : \c(SO_{-1}(3)) & \longrightarrow M_2(\C) \otimes M_2(\C) \otimes \c(SO(3)) \\ 
 a_{ij} & \longmapsto \tau_i \otimes \tau_j \otimes x_{ij}
 \end{align*}
in Theorem 5.2. Let $\tau_0$ be the identity matrix. Then 
$(\tau_i \otimes \tau_j)_{0 \leq i,j \leq 3}$ is a basis of $M_2(\C) \otimes M_2(\C)$.
Hence for $x \in H'$, there exists elements $p_{ij} \in H$ such that
$\theta (x) = \sum_{i,j} \tau_i \otimes \tau_j \otimes p_{ij}$. Now assume that $x \in J$. For $k \not = l$,
we have $a_{kl}x \in J \cap I_d'$, hence $a_{kl}x = 0$, and
$$\theta(a_{kl} x) = \sum_{i,j} \tau_k \tau_i \otimes \tau_l \tau_j \otimes x_{kl}p_{ij} =0$$
It follows that $x_{kl}p_{ij} =0$ for all $i$ and $j$.
The algebra $H$ is an integral domain (the algebraic 
group $SO(3)$ is irreducible), hence $p_{ij} =0$, $\forall i,j$. Therefore 
$\theta(x)=0$ and $x=0$.  
\end{proof}

Therefore, to complete the classification, it remains to find the quantum subgroups 
of the following quantum groups.

\begin{enumerate}
 \item $O_{-1}(2)$.
 \item $D_n^\tau$, $n$ even and $n>4$.
 \item $S_4^\tau$, the only non-trivial twist of the symmetric group $S_4$.
 \item $A_5^\tau$, the only non-trivial twist of the alternating group $A_5$.
\end{enumerate}

For the quantum orthogonal group $O_{-1}(2)$, this will be done in the next section.
The twisted dihedral groups are quantum subgroups of $O_{-1}(2)$,
and hence the classification of their quantum subgroups will follow from
the one for $O_{-1}(2)$.

Thus it remains to examine the twists $S_4^\tau$ and $A_5^\tau$. We use in a crucial way some results of Etingof and Gelaki \cite{eg2} to prove that their proper
quantum subgroups are all classical groups. 

\begin{lemma}
 The proper quantum subgroups of $S_4^\tau$ are all classical subgroups of $D_4$.
\end{lemma}

\begin{proof}
 We determine the algebra structure of the semisimple algebra $\c(S_4^\tau)$.
 First we know from \cite{bi0} that the maximal abelian quotient of 
  $\c(S_4^\tau)$ is $\c(D_4)$ and in particlar $C(S_4^\tau)$ has 8 non-equivalent
  one-dimensional representations.
  The 2-cocycle is induced by the Klein subgroup
  $$\V=\{(1,(1,2),(3,4),(1,2)(3,4)\}$$
  and is in minimal over $\c(\V)$ in the sense of \cite{eg2}.
  Let $g=(1,2,3)$. Following \cite{eg2}, we put $\V_g = \V \cap g\V g^{-1}$.
  We have $\V_g= \{1\}$ and hence Theorem 3.2 in \cite{eg2} furnishes
  an irreducible representation of $\c(S_4^\tau)$ having a dimension that is a
  multiple of $4$. An immediate counting argument (or Proposition 3.1)
  shows that this representation has dimension 4. Therefore we have an algebra isomorphism 
  $$\c(S_4^\tau) \simeq \C^8 \oplus M_4(\C)$$ 
  Now if $L$ is a non-commutative Hopf algebra quotient of $\c(S_4^\tau)$, its dimension
  is at least 16 and divides 24 by the Nichols-Zoeller theorem \cite{nz}: we are done.
\end{proof}

\begin{lemma}
The proper quantum subgroups of $A_5^\tau$ are all classical subgroups of $A_4$. 
\end{lemma}

\begin{proof} Again we begin by determining the algebra structure of
 $\c(A_5^\tau)$. Consider the Klein subgroup of $A_5$.
 $$\V=\{1, (1,2)(3,4), (1,3)(2,4), (1,4)(2,3)\}$$
The group inclusions $V \subset A_4 \subset A_5$ yield surjective
Hopf algebra maps $\c(A_5) \to \c(A_4) \to \c(\V)$. Twisting by the unique non-trivial 
2-cocycle on $V$ and using Lemma 6.2 and Lemma 6.4 now yields surjective Hopf algebra maps
$$\c(A_5^\tau) \to \c(A_4) \to \c(\V)$$
 This shows that $\c(A_5^\tau)$ has 12 non-equivalent one-dimensional representations
 and that the classical subgroups of $A_5^\tau$ are 
 the subgroups of $A_4$.
 Let $g_1 = (1,2,5)$ and put $\V_{g_1}= \V \cap g_1\V g_1^{-1}$.
 We have $V_{g_1}=\{1\}$ and hence Theorem 3.2 in \cite{eg2}
 furnishes an irreducible representation $V_1$ of $\c(A_5^\tau)$
 having a dimension that is a multiple of 4, and since $(\dim V_1)^2<60$, we have 
 $\dim V_1 = 4$. Now consider $g_2=(1,3,5)$ and $g_3=(1,4,5)$.
 The same reasonning as above furnishes irreducible 4-dimensional representations
 $V_2$ and $V_3$. For $i \not = j$ we have $\V g_i\V \not =\V g_j\V$ and by
 using \cite{eg2}, we see that the representations $V_1$, $V_2$ and $V_3$ are pairwise non-equivalent. We conclude that
 $$\c(A_5^\tau) \simeq \C^{12} \oplus M_4(\C)^3$$ 
 Let $L$ be a non-commutative Hopf algebra quotient of $\c(A_5^\tau)$:
 its dimension is at least 16 and divide 60 by the Nichols-Zoeller Theorem,
 thus $\dim L =20,30$ or $60$. A counting argument shows that $\dim L \neq 30$. 
 Assume that $\dim L = 20$. Then, as an algebra, $L \simeq \C^4\oplus M_4(\C)$.
 Hence the maximal abelian quotient of $L$ is $\C^4 \simeq \c(\V)$, and we have Hopf algebra surjections $\c(A_5^\tau) \to L \to \c(\V)$. Twisting
 by the appropriate 2-cocycle yields Hopf algebra surjections
 $\c(A_5) \to L' \to \c(\V)$. This means that $A_5$ has a subgroup of order 20, which is known to be false. Hence $\dim L=60$ and we are done. 
\end{proof}

\section{The quantum group $O_{-1}(2)$ and its quantum subgroups}

In this section we classify the compact quantum subgroups of $O_{-1}(2)$, finishing in this way the proof of the classification theorem. The novelty
will be the occurence of $\widehat{D}_\infty$, the quantum dual of the infinite
dihedral group
and of pseudo-twisted dicyclic groups $DC_n^\tau$.

First let us note that we can classify the  
quantum subgroups diagonally containing the Klein subgroup
in a manner similar to the one of the previous section. In particular
we conclude that the quantum groups $D_n^\tau$, $n>4$ and even, all
are quantum subgroups of $O_{-1}(2)$. But the techniques of the end of the last
section do not work here, because $O(2)$ is not connected.
So we will proceed directly, and another presentation of $\r(O_{-1}(2))$
will be more convenient.

Recall \cite{bbc} that $A_h(2)$ is the universal $*$-algebra
 generated by the entries of a cubic matrix. More precisely $A_h(2)$
 is the universal $*$-algebra
 presented by generators
$(v_{ij})_{1\leq i,j\leq n}$ and relations:
\begin{enumerate}
 \item The matrix  $v=(v_{ij})$ is orthogonal.
 \item $v_{ij}v_{ik}= v_{ik}v_{kj} = 0=v_{ji}v_{ki} = v_{ki}v_{ji}$ if $j \not = k$. 
\end{enumerate}
The Hopf $*$-algebras $\r(O_{-1}(2))$ and $A_h(2)$ are isomorphic, via
$$
\left(
\begin{array}{ll}
a_{11} & a_{12} \\ 
a_{21} & a_{22}
\end{array}\right) \longmapsto 
\left(
\begin{array}{ll}
\frac{1}{\sqrt{2}} & \frac{-1}{\sqrt{2}}\\ 
\frac{1}{\sqrt{2}} & \frac{1}{\sqrt{2}}
\end{array}\right)
\left(
\begin{array}{ll}
v_{11} & v_{12} \\ 
v_{21} & v_{22}
\end{array}\right)
\left(
\begin{array}{ll}
\frac{1}{\sqrt{2}} & \frac{1}{\sqrt{2}} \\ 
\frac{-1}{\sqrt{2}} & \frac{1}{\sqrt{2}}
\end{array}\right)$$
It is shown in \cite{bi2} that there is an algebra isomorphism
$$A_h(2) \simeq \C[\z_2*\z_2] \otimes \c[\z_2]$$
Using the group isomorphism $\z_2 * \z_2 \simeq D_\infty$,
we get the following isomorphism, that will enable us
to connect $A_h(2)$ with the algebras constructed by Masuoka in \cite{ma}. 

\begin{lemma}
 Let $Z$ be the universal algebra presented by generators
 $z,z^{-1},t,a$ and relations
 $$zz^{-1}=1 = z^{-1}z, \ t^2=1=a^2, \ tz = z^{-1}t, \ at=ta, \ az= za$$
 We have an algebra isomorphism
  \begin{align*}
  A_h(2) & \longrightarrow Z\\
  \left(
\begin{array}{ll}
v_{11} & v_{12} \\ 
v_{21} & v_{22}
\end{array}\right)
& \longmapsto
\left(
\begin{array}{ll}
t\frac{(1+a)}{2} & t\frac{(1-a)}{2} \\ 
tz\frac{(1-a)}{2} & tz\frac{(1+a)}{2}
\end{array}\right)
\end{align*}
\end{lemma}

\begin{proof}
 The map below defines the inverse isomorphism.
 $$t \mapsto v_{11} +v_{12}, \ z \mapsto v_{11}v_{22} + v_{12}v_{21}, \ a 
 \mapsto v_{11}^2 -v_{12}^2$$ 
\end{proof}

We now describe the Peter-Weyl decomposition of $A_h(2)$. We put $d=v_{11}^2-v_{12}^2$. This a group-like element satisfying $d^2=1$. For $k\in\mathbb N^*$, we define a comatrix subcoalgebra of $A_h(2)$ in the following way. For $k=2m$ we put
$$\begin{pmatrix}v_{11}^{[k]}&v_{12}^{[k]}\cr v_{21}^{[k]}&v_{22}^{[k]}\end{pmatrix}=\begin{pmatrix}(v_{11}v_{22})^m&(v_{12}v_{21})^m\cr (v_{21}v_{12})^m&(v_{22}v_{11})^m\end{pmatrix}$$
For $k=2m+1$, we put
$$\begin{pmatrix}v_{11}^{[k]}&v_{12}^{[k]}\cr v_{21}^{[k]}&v_{22}^{[k]}\end{pmatrix}=\begin{pmatrix}(v_{11}v_{22})^mv_{11}&(v_{12}v_{21})^mv_{12}\cr (v_{21}v_{12})^mv_{21}&(v_{22}v_{11})^mv_{22}\end{pmatrix}$$
We put $C(k)=\mathbb Cv_{11}^{[k]}\oplus\mathbb Cv_{12}^{[k]}\oplus\mathbb Cv_{21}^{[k]}\oplus\mathbb Cv_{22}^{[k]}$. The vector space $C(k)$ has dimension 4 by Lemma 7.1, and is a (simple) subcoalgebra of $A_h(2)$. An application of Lemma 7.1 yields the Peter-Weyl decomposition of $A_h(2)$.

\begin{proposition}
The Peter-Weyl decomposition of $A_h(2)$ is
$$A_h(2)=\C 1\oplus\C d\oplus(\oplus_{k\geq 1}C(k))$$ 
\end{proposition}

We now construct some finite-dimensional quotients of $A_h(2)$.
For $m \in \mathbb N^*$ and $e= \pm 1$ we define $A(2m,{e})$ to be the quotient
of $A_h(2)$ by the relations
$$(v_{11}v_{22})^m = (v_{22}v_{11})^m, \ 
(v_{12}v_{21})^m = e(v_{21}v_{12})^m$$
Similarly, for $m \in \mathbb N$ and $e=\pm 1$,
we define $A(2m+1,{e})$ to be the quotient
of $A_h(2)$ by the relations
$$(v_{11}v_{22})^mv_{11} = (v_{22}v_{11})^mv_{22}, \ 
(v_{12}v_{21})^m v_{12}= e(v_{21}v_{12})^mv_{21}$$

\begin{proposition}
 The algebras $A(k, {e})$ are Hopf $*$-algebra quotients of $A_h(2)$ for any $k \in \mathbb N^*$.
 \begin{enumerate}
  \item For $k=1$, we have $A(1,{e})\simeq\C[\V]$.
  \item For $k=2$, the Hopf algebra $A(2,1)$ has dimension $8$ and is isomorphic
  with $\c(D_4)$.
  \item For $k>2$, the Hopf algebra $A(k,1)$ has dimension $4k$ and is isomorphic
  with $\c(D_{2k}^\tau)$, the function algebra on $D_{2k}^\tau$, the unique 
 non-trivial twist of the dihedral group $D_{2k}$.  
  \item For $k>2$, the Hopf algebra $A(k,-1)$ has dimension $4k$ and is isomorphic
  with $\c(DC_{k}^\tau)$, the function algebra on $DC_{k}^\tau$, a non-trivial 
 pseudotwist of the dicyclic group $DC_{k}$ considered by Masuoka \cite{ma},
 Nikshych \cite{ni}, Suzuki \cite{su} and Vainerman \cite{va}.  
 \end{enumerate}
\end{proposition}

\begin{proof}
 It is straightforward to check that $A(k,{e})$ is indeed a Hopf $*$-algebra quotient of $A_h(2)$, and the assertion at $k=1$ is immediate since $A(1, {e})$ is
 commutative. Also $A(2,1)$ is commutative, and hence is the function algebra on $D_4$.
 
 Assume that $k>2$. Then by using Lemma 7.1, we see that $A(k,1)$
 is isomorphic with the quotient of the algebra $Z$ of Lemma 7.1
 by the relation $z^k=1$. Therefore $A(k,1)$ is isomorphic with the algebra
 $\mathcal A_{4k}$ considered by Masuoka (Definition 3.3 in \cite{ma}).
 This is a Hopf algebra isomorphism, and we conclude by Theorem 4.1 in \cite{ma}
 that $A(k,1) \simeq \c(D_{2k}^\tau)$.
 
 Similarly, by  using Lemma 7.1, we see that $A(k,-1)$
 is isomorphic with the quotient of the algebra $Z$ 
 by the relation $z^k=a$. Therefore $A(k,-1)$ is isomorphic with the algebra
 $\mathcal B_{4k}$ considered by Masuoka (Definition 3.3 in \cite{ma}).
 This is a Hopf algebra isomorphism. In \cite{ni,va} the Hopf algebra
 $\mathcal B_{4k}$ appears as a pseudo-twist of the dicyclic group,
 and this is why we put $\mathcal B_{4k} = \c(DC_{k}^\tau)$.
\end{proof}

We also will need the following technical lemma.

\begin{lemma}
 Let $A$ be a Hopf $*$-algebra and let $w=(w_{ij})$ in $M_2(A)$ be a unitary matrix
 such 
that $\Delta(w_{ij})=\sum_kw_{ik}\otimes w_{kj}$ and $\varepsilon(w_{ij})=\delta_{ij}$. Assume
that 
$w_{ij}w_{ik}=w_{ik}w_{ij}=0=w_{ji}w_{ki}=w_{ki}w_{ji}$ if $j\neq k$ and that
there exists $F\in GL(2)$ such that the matrix $FwF^{-1}$ is diagonal. 

\noindent
{\rm (A)} Assume that the matrix $w$ is orthogonal. Then one of the
following happens.
\begin{enumerate}
\item $w_{12}=0=w_{21}$ and $w_{11}^2=1=w_{22}^2$.
\item $w_{11}=w_{22}$ and $w_{12}=\pm w_{21}$.
\end{enumerate}
If moreover 
$$FwF^{-1} = \
\begin{pmatrix}
 1 & 0\\
 0 & *
\end{pmatrix}$$
Then according to the previous cases, one the following happens.
\begin{enumerate}
 \item $w_{11}=1$ or $w_{22}=1$.
 \item $w_{11}+w_{22} = 1$ or $w_{11}+iw_{22}=1$
\end{enumerate}
\noindent
{\rm (B)} Assume that the coefficients of the matrix $w$ pairwise commute
and that $w_{11}w_{22} + w_{12}w_{21} =1$. Then one of the following happens.
\begin{enumerate}
\item $w_{12}=0=w_{21}$ and $w_{11}w_{22}=1=w_{22}w_{11}$.
\item $w_{11}=w_{22}$ and $w_{12}=\lambda w_{21}$ for $\lambda\in\mathbb C$ with $|\lambda|=1$.
\end{enumerate}

\end{lemma}

\begin{proof}
Let $F=(^a_c{\ }^b_d)$. Then if $FwF^{-1}$ is diagonal, we have the following equations
$$ab(w_{11}-w_{22})=a^2w_{12}-b^2w_{21}$$
$$cd(w_{11}-w_{22})=c^2w_{12}-d^2w_{21}$$
Assume that we are in situation (A).
If $b=0$, we have $a^2w_{12}=0$ and hence $w_{12}=0$. By applying the antipode we get $w_{21}=0$, and finally the orthogonality condition ensures that we are in case (1). The same reasoning holds if $a=0$, $c=0$ or $d=0$.

Assume now that the scalars $a,b,c,d$ are all nonzero. Then by multiplying the first equation in the proof by $w_{11}$ we get
$$w_{11}^2=w_{11}w_{22}=w_{22}w_{11}=w_{22}^2$$
and then
$$w_{11}=w_{11}^3=w_{11}w_{22}w_{22}=w_{22}^3=w_{22}.$$
We also get
$$(ab^{-1}-cd^{-1})w_{12}=(ba^{-1}-dc^{-1})w_{21}$$
and hence there exists $\lambda\in\mathbb C^*$ such that $w_{12}=\lambda w_{21}$. We then use the antipode to get $\lambda=\pm 1$, and we are in case (2).

In case (1), the group-like elements in the coalgebra generated by the entries
of $w$ are $w_{11}$ and $w_{22}$, while in case (2) they are
$w_{11}+w_{22}$ and $w_{11}+iw_{22}$. The last assertion follows.

Let us assume now that condition (B) holds.
If $b=0$, then $a^2w_{12}=0$ and $w_{12}=0=w_{21}$: we are in case (1), and the same reasoning holds if $a$, $c$ or $d$ is zero.

Assume now that the scalars $a,b,c,d$ are all non-zero. Similarly to the previous case, we get $w_{12}=\lambda w_{21}$ for some $\lambda\in\mathbb C^*$, and by using the adjoint we get $|\lambda|=1$. Finally, since $S(w_{11})=w_{22}$ and $S(w_{12})=w_{12}$, by applying the antipode to the first equation at the beginning of the proof yields $w_{11}=w_{22}$: we are in situation (2).
\end{proof}

We now turn to classify the Hopf $*$-algebra quotients of $A_h(2)$.
Let $L$ be a non-trivial Hopf $*$-algebra quotient of $A_h(2)$. For simplicity, the projections of the elements of $A_h(2)$ are denoted by the same symbol in $L$.
The Peter-Weyl decomposition of $A_h(2)$ becomes 
a coalgebra sum
$$L = \C 1 + \C d + (\sum_{k\geq 1}C(k))$$ 
For $k \in \mathbb N^*$ the $2$-dimensional corepresentation
corresponding to the coalgebra $C(k)$ is denoted $V_k$, and $V_k$ is irreducible
if and only if $\dim C(k)=4$.

\begin{lemma}
 Assume that $d=v_{11}^2-v_{12}^2=1$ in $L$.
 Then $L$ is isomorphic with $\C[D_\infty]$, the group
 algebra of the infinite dihedral group or with
 $\C[D_n]$, the group algebra of the dihedral group of order
 $2n$, for $n \in \mathbb N^*$.
\end{lemma}

\begin{proof}
 Summing $d$ with $v_{11}^2+v_{12}^2=1$, we get $v_{11}^2=v_{22}^2=1$ and $v_{12}^2=v_{21}^2=0$. Since $v_{12}$ and $v_{21}$ are self-adjoint, we have $v_{12}=0=v_{21}$. We get a surjective Hopf $*$-algebra map
 $\C[D_\infty] \simeq \C[\z_2*\z_2] \to L$, and we conclude 
 by using the classification of Hopf $*$-algebra quotients
 of $\C[D_\infty]$.
\end{proof}

So we assume now that $d \not = 1$ in $L$, and that $L$ is not a quotient of 
$\C[D_\infty]$.

\begin{lemma}
 Assume that the corepresentation $V_k$ is not irreducible
 for some $k \in \mathbb N^*$.
 Then there exists a surjective Hopf $*$-algebra map
 $A(k,e) \to L$ for $e=\pm 1$, and in particular
 $\dim L$ divides $4k$.  
\end{lemma}

\begin{proof}
 Assume first that $k=2m+1$ is odd.
Consider the matrix 
 $$\begin{pmatrix}(v_{11}v_{22})^mv_{11}&(v_{12}v_{21})^mv_{12}\cr (v_{21}v_{12})^mv_{21}&(v_{22}v_{11})^mv_{22}\end{pmatrix}$$
 Then by the first part of Lemma 7.2 one of the following holds:
 \begin{enumerate}
 \item $(v_{12}v_{21})^mv_{12} = 0 = (v_{21}v_{12})^mv_{21}$ and
 $((v_{11}v_{22})^mv_{11})^2 = 1 = ((v_{22}v_{11})^mv_{22})^2$.
 \item  $(v_{11}v_{22})^mv_{11} = (v_{22}v_{11})^mv_{22}$ and
 $(v_{12}v_{21})^mv_{12} = \pm (v_{12}v_{21})^mv_{12}$.
\end{enumerate}
In the first situation some direct computations
give $v_{11}^2=1=v_{22}^2$ and $v_{12} = v_{21}=0$, contradicting the assumption on $L$.
The second family of  relations 
yields a surjective Hopf $*$-algebra map $A(k,e) \to L$.

Similarly, if $k=2m$ is even, we consider the matrix
$$\begin{pmatrix}(v_{11}v_{22})^m&(v_{12}v_{21})^m\cr (v_{21}v_{12})^m&(v_{22}v_{11})^m\end{pmatrix}$$
By the second part of Lemma 7.2, one of the following holds:
 \begin{enumerate}
 \item $(v_{12}v_{21})^m = 0 = (v_{21}v_{12})^m$ and
 $((v_{11}v_{22})^m)^2 = 1 = ((v_{22}v_{11})^m)^2$.
 \item  $(v_{11}v_{22})^m = (v_{22}v_{11})^m$ and
 $(v_{12}v_{21})^m = \lambda (v_{21}v_{12})^m$.
\end{enumerate}
Similarly to the previous case, the first family 
of relations cannot hold. If $\lambda = \pm 1$ in the second family of relations,
we have our result. Otherwise we have
\begin{align*}
(v_{12}v_{21})^mv_{21} & = \lambda (v_{21}v_{12})^mv_{21}\\
\Rightarrow  (v_{12}v_{21})^{m-1} v_{12} &= \lambda v_{21}(v_{12}v_{21})^m
= \lambda^2 v_{21}(v_{21}v_{12})^{m}\\
&= \lambda^2v_{21}^2(v_{12}v_{21})^{m-1}v_{12} = 
 \lambda^2 (v_{12} v_{21})^{m-1}v_{12}
\end{align*} 
If $\lambda^2 \not = 1$, we have then  
$$(v_{12} v_{21})^{m-1}v_{12} =0=(v_{12} v_{21})^m= (v_{21}v_{12})^m$$
and hence $(v_{12} v_{21})^m= \pm(v_{21}v_{12})^m$: we have our result.

The last assertion follows from the Nichols-Zoeller theorem. 
\end{proof}

This last lemma before the classification of the
Hopf $*$-algebra quotients of $A_h(2)$ uses the following fusion rules of 
$A_h(2)$ (see \cite{bi2}, these are the same as the ones of $O(2)$):
$$V_i \otimes V_j \simeq V_{i+j} \oplus V_{|i-j|} \ if
\ i \not = j \in \mathbb N^* \ ;$$
$$V_i \otimes V_{i} \simeq\C \oplus \C d \oplus  
V_{2i} \ ; \  \C d\otimes  \C d\simeq  \C \ ; \
V_i \otimes  \C d\simeq  \C d\otimes V_i \simeq V_i, \ i \in \mathbb N^*.$$ 

\begin{lemma}
 Let $j >2$ and assume that
 \begin{enumerate}
  \item The corepresentations $V_1, \ldots , V_{j-1}$ of $L$ are
  irreducible and pairwise non-equivalent.
  \item For $l \in \{1, \ldots, j-1\}$, we have
  $V_l\simeq V_j$.
 \end{enumerate}
Then $j \not = l+1$ and $L \simeq A(l+j,e)$ for $e=\pm 1$.
\end{lemma}

\begin{proof}
 We have $V_l \simeq V_j$ for $l<j$ and hence
 $$V_l \otimes V_j \simeq V_{j-l} \oplus V_{l+j}
 \simeq V_l \otimes V_l \simeq \C \oplus \C d \oplus V_{2l}$$
 Hence $V_{l+j}$ is not irreducible, with $V_{l+j} \simeq \C \oplus \C d$, while
 $V_{2l}$ is irreducible with $V_{2l} \simeq V_{j-l}$.
 By the previous lemma we have a surjective
 Hopf $*$-algebra map $A(l+j,e) \to L$ and $\dim L$ divides $4(l+j)$.
 A counting argument using the first assumption shows that 
 $\dim L \geq 4(j-1)+2 = 4j-2$. 
 
 Assume that $j-l>1$. Then $4j-2>4(l+j)/2$, and hence 
 $\dim L > 4(l+j)/2$. We conclude that $\dim L = 4(l+j)$ and that 
 $L \simeq A(l+j,e)$.
 
 Assume now that $j=l+1$. Then $V_{2l+1} \simeq \C \oplus \C d$ and $V_{2l}$ is irreducible. By Lemma 7.2, we have
 $(v_{11}v_{22})^lv_{11} = 1$, or   $(v_{22}v_{11})^lv_{11} = 1$, or
 $(v_{11}v_{22})^lv_{11} +(v_{22}v_{11})^lv_{22} = 1$
 or $(v_{11}v_{22})^lv_{11} +i(v_{22}v_{11})^lv_{22} = 1$.
 If for example $(v_{11}v_{22})^lv_{11} = 1$, we get, multiplying by $v_{11}$,
 that $(v_{11}v_{22})^l= (v_{22}v_{11})^l =v_{11}$, which contradicts the irreducibility
 of $V_{2l}$. The other cases are treated similarly, and we are done.
\end{proof}

We are now ready to complete the classification of the compact quantum subgroups of $O_{-1}(2)$.

\begin{theorem} 
Any non-trivial Hopf $*$-algebra quotient of $A_h(2)$ is isomorphic to one of the following Hopf $*$-algebras.
\begin{enumerate}
 \item $\C[D_\infty]$, the group algebra of the infinite dihedral group.
 \item $\C[D_{n}]$, the group algebra of the dihedral group of order $2n$, $n\in \mathbb N^*$.
 \item $A(k,{e})$, $k\in\mathbb N^*$, ${e}=\pm 1$.
\end{enumerate}
\end{theorem}

\begin{proof}
Let $L$ be a non-trivial Hopf $*$-algebra quotient of $A_h(2)$.
We can assume that $L$ is not a quotient of $\C[D_\infty]$, and that $d\not =1$.
Then the comodule $V_1$ is irreducible. 
If all the comodules $V_i$ are irreducible and pairwise non-equivalent, 
it is clear that $L \simeq A_h(2)$.
Hence we can consider a family of pairwise non-equivalent irreducible comodules
$V_1, \ldots , V_{j-1}$ ($j >1$) such that either
\begin{enumerate}
 \item $V_j$ is irreducible and there exists $l \in \{1,\ldots , j-1\}$
 with $V_l \simeq V_j$,
 \item $V_j$ is not irreducible. 
\end{enumerate}
The first situation is the one of the previous lemma, and we
conclude that $L \simeq A(l+j,e)$. 

In the second situation, we have, by Lemma 7.4, 
a surjective Hopf $*$-algebra map $A(j,e) \to L$, and $\dim L | 4j$.
But by the assumption $\dim L \geq 4j-2$, and hence
we conclude that $\dim L = 4j$ and that $L \simeq A(j,e)$.  
\end{proof}

\section{Strategy for the ADE classification: algebraic invariants}

The proof of the main result, Theorem 1.1, is now complete. In the rest of this paper we 
describe an ADE type classification of the quantum subgroups we have found.
This section is devoted to the description of our strategy.

The ADE classification, i.e. McKay correspondence, for (compact) subgroups
of $SO(3)$ is obtained as follows. Starting from such a subgroup $G \subset SO(3)$,
we consider its double cover $G'\subset SU(2)$, and associate to it its 
McKay graph \cite{mc} to obtain an ADE graph. Here $\Q_4 \simeq SO_{-1}(3)$
does not have such a double cover, and we do not know a direct 
way to associate a graph to a subgroup of $\Q_4$. Therefore we proceed as follows. 
 
We consider the following two types of objects.

\begin{enumerate}
\item Quantum groups acting on 4 points.
\item ADE type graphs.
 \end{enumerate}

To any such object $X$ (a quantum group or a graph), we associate
a sequence of integers $(c_n)$, or equivalently a probability measure
on the circle (the circular measure), that we call the algebraic
invariants of  $X$.
Then we say that a quantum group and a graph correspond to each other
if they have the same algebraic invariants.

This might look a bit artificial, but the following facts
support our strategy.

\begin{enumerate}
\item The classical McKay correspondence might be described in such a way, see
Section 11.
\item When the quantum group acts ergodically on $\C^4$, one can associate a subfactor
of index $4$ 
to this data \cite{ba1}. The principal graph of the subfactor
is the ADE graph we find.
\end{enumerate}

Let us begin with the algebraic invariants of a quantum permutation group 
$\mathcal G\subset\mathcal Q_4$. These depend not only on $\mathcal G$, but also on its embedding into $\mathcal Q_4$. These invariants come from the spectral theory of the character of the fundamental representation of $\mathcal G$, on the space $\mathbb C^4$. They are best introduced in the following way.

\begin{definition}
The algebraic invariants of a quantum subgroup $\mathcal G\subset\mathcal Q_4$ are:
\begin{enumerate}
\item The multiplicity $c_k$. This is the number of copies of $1$ into $u^{\otimes k}$.
\item The spectral measure $\mu$. This is the measure having $c_k$ as moments.
\item The circular measure $\varepsilon$. This is the pullback of $\mu$ via $\varphi(q)=(q+\bar{q})^2$.
\end{enumerate}
\end{definition}

In this definition $\mu$, $\varepsilon$ are probability measures on the real line, and on the unit circle. The support of $\mu$ being contained in $[0,4]$ and the map $\varphi$ being surjective from the unit circle to $[0,4]$, the measure $\varepsilon$ is well-defined, by distributing the weight at $x\in [0,4]$ uniformly among its preimages $q=\varphi^{-1}(x)$. See \cite{ba2}, \cite{bbi}.

We use the following notations:
\begin{enumerate}
\item $d$ is the uniform measure on the unit circle. \item $d_n$ is the uniform measure on $2n$-roots of 1.
\item $d_n'$ is the uniform measure on odd $4n$-roots of 1.
\item $e_n=d_1'$ if $n$ is even and $e_n=d_2'$ if $n$ is odd.
\item For $s \in [0,4]$, $\gamma_s$ is the pullback of $\delta_s$, the Dirac mass
at $s$, via $\varphi(q) = (q + \bar{q})^2$.
\end{enumerate}

We let $\alpha(q)=2Im(q)^2$. This is the pullback via $\varphi$ of the free Poisson density  $(2\pi)^{-1}\sqrt{4x^{-1}-1}$ of Marchenko-Pastur \cite{mp} and Voiculescu \cite{vdn}. 

\smallskip

We discuss now algebraic invariants of graphs. 
The interest in these invariants might seem less clear than in the quantum group case. The idea is that these are related to the quantum group or subfactor ones via the representation-theoretic notion of principal graph. See \cite{bbi,jo2}. 

\begin{definition}
The algebraic invariants of a rooted graph of norm $\leq 2$ are:
\begin{enumerate}
\item The loop number $c_k$. This is the number of $2k$-loops based at the root.
\item The spectral measure $\mu$. This is the measure having $c_k$ as moments.
\item The circular measure $\varepsilon$. This is the pullback of $\mu$ via $\varphi(q)=(q+\bar{q})^2$.
\end{enumerate}
\end{definition}

The comment after previous definition applies as well to this situation, and shows that $\varepsilon$ is indeed a probability measure on the unit circle. Here the fact that $\mu$ is supported by $[0,4]$ comes from the fact that the norm of the graph is $\leq 2$. 

In fact, $\mu$ is the spectral measure of the square of the adjacency matrix of the graph, and the norm $\leq 2$ condition means that this adjacency matrix is supported by $[-2,2]$, hence its square is supported by $[0,4]$. See \cite{bbi} for details.

\medskip

The remaining sections are devoted to the description of the circular
measures of quantum permutation groups and ADE graphs.
For this purpose 
it is convenient to divide the quantum subgroups of $\Q_4$ in several classes,
which leads to the 
following  reformulation of Theorem 1.1.

\begin{theorem}
The quantum subgroups of $\mathcal Q_4$ are of three types.
\begin{enumerate}
\item Classical: $\mathbb Z_1$, $\mathbb Z_2$, $\mathbb Z_3$, $\mathbb Z_4$, $V$, $S_3$, $D_4$, $A_4$, $S_4$.
\item Quantum, single: $S_4^\tau$, $A_5^\tau$, $\widehat{D}_\infty$, $O_{-1}(2)$, $\mathcal Q_4$.
\item Quantum, series: $\widehat{D}_{n}$, $D_{2n}^\tau$ with $n\geq 3$, $DC_n^\tau$ with $n\geq 2$.
\end{enumerate}
\end{theorem}

We begin by treating the case of the subgroups of $S_4$ in the next section.
In Section 10, we describe the circular measures of the extended Coxeter-Dynkin
graphs (taken from \cite{bbi}), plus two other graphs. In Section 11, we treat
the remaining quantum subgroups, using a reformulation of McKay correspondence
at the level of algebraic invariants, and Section 12 is devoted to the final classification tables.

\section{Algebraic invariants of subgroups of $S_4$}

The subgroups of $S_4$ are $\mathbb Z_1$, $\mathbb Z_2$, $\mathbb Z_3$, $\mathbb Z_4$, $V$, $S_3$, $D_4$, $A_4$, $S_4$.
In this section we discuss the algebraic invariants of these groups. 

Before starting, we have to discuss a quite undexpected issue.

\begin{definition}
The $2$-element groups acting on $4$ points are as follows.
\begin{enumerate}
\item $D_1$: here the non-trivial element acts with $2$ fixed points.
\item $\mathbb Z_2$: here the non-trivial element acts without fixed point.
\end{enumerate}
The Klein groups acting on $4$ points are as follows.
\begin{enumerate}
\item $D_2$: here the non-trivial elements acts without fixed points.
\item $V$: here there are $2$ non-trivial elements having $2$ fixed points.
\end{enumerate}
\end{definition}

The group $D_1$ is to be added to the above list, as a $10$-th classical subgroup of $S_4$. Observe that the notation $D_1$ is in agreement with the isomorphisms for small dihedral groups, which are $D_1\simeq\mathbb Z_2$, $D_2\simeq V$, $D_3\simeq S_3$.

The point is that the algebraic invariants of a quantum permutation group $\mathcal G\subset\mathcal Q_4$ depend not only on $\mathcal G$, but also on its embedding into $\mathcal Q_4$. 
In fact $D_2$ has the same algebraic invariants as $\mathbb Z_4$,
so we will ignore it.

The computation for symmetric groups leads to the following formulae, the notations
for circular measures being the ones of the previous section.

\begin{theorem}
The circular measures of subgroups $G\subset S_4$ are as follows.
\begin{enumerate}
\item For $\mathbb Z_1$, $\mathbb Z_2$, $\mathbb Z_3$, $V$ we have $\varepsilon=d_n$ with $n=1,2,3,4$.
\item For $D_1$, $\mathbb Z_4$, $S_3$, $D_4$ we have $\varepsilon=(e_n+d_n)/2$ with $n=1,2,3,4$.
\item For $A_4$, $S_4$ we have  $\varepsilon=\alpha\,d_n+(d_{n-1}-d_n)/2$ with $n=3,4$.
\end{enumerate}
\end{theorem}

\begin{proof}
For $G\subset S_4$ we denote by $m_s$ the number of elements in $G$ having $s$ fixed points. Let also $m=|G|$. Then the invariants of $G$ are:
\begin{eqnarray*}
c_k&=&\frac{1}{m}\sum_{s=0}^4m_s\,s^k\cr
\mu&=&\frac{1}{m}\sum_{s=0}^4m_s\,\delta_s\cr
\varepsilon&=&\frac{1}{m}\sum_{s=0}^4m_s\,\gamma_s
\end{eqnarray*}

These formulae are all equivalent, and the first two ones are valid in fact for any subgroup $G\subset S_n$, with $4$ replaced of course by $n$. See \cite{ba2} for details.

We have $m_3=0$ and $m_4=1$, so all invariants can be expressed in terms of the numbers $m_0,m_1,m_2$ and $m=m_0+m_1+m_2+1$:
\begin{eqnarray*}
c_k&=&\frac{1}{m}(m_1+m_2\,2^k+4^k)\cr
\mu&=&\frac{1}{m}(m_0\,\delta_0+m_1\,\delta_1+m_2\,\delta_2+\delta_4)&\cr
\varepsilon&=&\frac{1}{m}(m_0\,\gamma_0+m_1\,\gamma_1+m_2\,\gamma_2+\gamma_4)
\end{eqnarray*}

The relevant measures $\gamma_s$ are as follows:
\begin{eqnarray*}
\gamma_0&=&d_1'\cr
\gamma_1&=&(3d_3-d_1)/2\cr
\gamma_2&=&d_2'\cr
\gamma_4&=&d_1
\end{eqnarray*}

Thus we have the following formula for the circular measure:
$$\varepsilon=\frac{1}{m}\left(m_0\,d_1'+\frac{m_1}{2}(3d_3-d_1)+m_2\,d_2'+d_1\right)$$

The numbers $m_0,m_1,m_2$ and $m=m_0+m_1+m_2+1$ are as follows:
\setlength{\extrarowheight}{3pt}
\begin{center}
\begin{tabular}[t]{|l|l|l|l|l|l|l|l|l|l|l|}
\hline &$\mathbb Z_1$&$\mathbb Z_2$&$\mathbb Z_3$&$V$&$D_1$&$\mathbb Z_4$&$S_3$&$D_4$&$A_4$&$S_4$\\ 
\hline $m_0$&0&1&0&1&0&3&0&5&3&9\\
\hline $m_1$&0&0&2&0&0&0&2&0&8&8\\
\hline $m_2$&0&0&0&2&1&0&3&2&0&6\\
\hline   $m$&1&2&3&4&2&4&6&8&12&24\\
\hline
\end{tabular}
\end{center}
\setlength{\extrarowheight}{0pt}
\bigskip

Observe that $V$ acts with $m_0=0$. The other action of $V$, with $m_0=3$, gives the same algebraic invariants as $\mathbb Z_4$, so we will not consider it.

We are now in position of computing the circular measures. We will use many times the formula $d_n+d_n'=2d_{2n}$. 

For the groups $\mathbb Z_1$, $\mathbb Z_2$, $\mathbb Z_3$, $V$, 
$D_1$, $\mathbb Z_4$, $S_3$, $D_4$ we immediatly get the measures in the statement.

For $A_4$ we have the following preliminary computation:
\begin{eqnarray*}
\varepsilon&=&(3d_1'+4(3d_3-d_1)+d_1)/12\cr
&=&(3d_1'-3d_1+12d_3)/12\cr
&=&(d_1'-d_1+4d_3)/4\cr
&=&d_3+(d_1'-d_1)/4
\end{eqnarray*}

For $S_4$ we have the following preliminary computation:
\begin{eqnarray*}
\varepsilon&=&(9d_1'+4(3d_3-d_1)+6d_2'+d_1)/24\cr
&=&(9d_1'-3d_1+6d_2'+12d_3)/24\cr
&=&(6d_2+6d_1'-6d_1+6d_2'+12d_3)/24\cr
&=&(12d_4+6d_1'-6d_1+12d_3)/24\cr
&=&(d_3+d_4)/2+(d_1'-d_1)/4
\end{eqnarray*}

Consider now the fundamental density $\alpha(q)=2Im(q)^2$. The measure $\alpha\,d_3$ is uniformly supported by the imaginary $6$-roots, so we have $\alpha\,d_3=\gamma_1$. We get:
\begin{eqnarray*}
\alpha\,d_3+(d_2-d_3)/2
&=&\gamma_1+(d_2-d_3)/3\cr
&=&(3d_3-d_1)/2+(d_2-d_3)/2\cr
&=&d_3+(d_2-d_1)/2\cr
&=&d_3+(d_1'-d_1)/4
\end{eqnarray*}

Also, regarding the measure $\alpha\,d_4$, this has density $1/4$ at $\pm i$, and density $1/8$ at the odd $8$-roots. Thus we have $\alpha\,d_4=(\gamma_0+\gamma_2)/2$, which gives:
\begin{eqnarray*}
\alpha\,d_4+(d_3-d_4)/2
&=&(\gamma_0+\gamma_2)/2+(d_3-d_4)/2\cr
&=&(d_1'+d_2')/2+(d_3-d_4)/2\cr
&=&(2d_2-d_1+2d_4-d_2)/2+(d_3-d_4)/2\cr
&=&(d_3+d_4)/2+(d_2-d_1)/2\cr
&=&(d_3+d_4)/2+(d_1'-d_1)/4\cr
\end{eqnarray*}

This gives the formulae for $A_4,S_4$ in the statement.
\end{proof}

\section{Algebraic invariants of ADE type graphs}

We describe now the algebraic invariants of ADE type graphs. 
The extended Coxeter-Dynkin diagrams are as follows:
$$\ \ \ \ \ \ \ \tilde{A}_n=
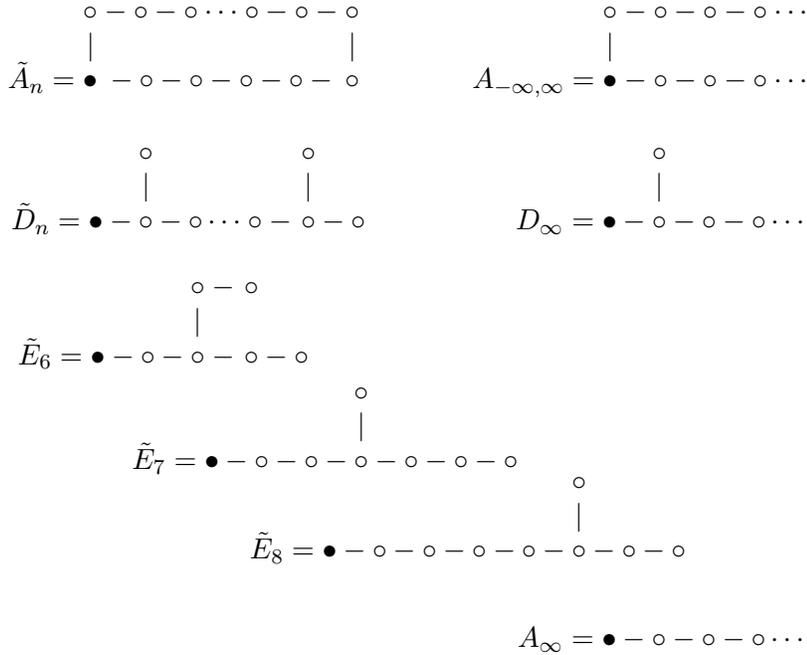
\begin{matrix}
\circ&\!\!\!\!-\circ-\circ\cdots\circ-\circ-&\!\!\!\!\circ\cr
|&&\!\!\!\!|\cr
\bullet&\!\!\!\!-\circ-\circ-\circ-\circ-&\!\!\!\!\circ\cr\cr\cr\end{matrix}\hskip15mm A_{-\infty,\infty}=
\begin{matrix}
\circ&\!\!\!\!-\circ-\circ-\circ\cdots\cr
|&\cr
\bullet&\!\!\!\!-\circ-\circ-\circ\cdots\cr\cr\cr\end{matrix}
\hskip15mm$$
\vskip-9mm
$$\;\tilde{D}_n=\bullet-
\begin{matrix}\circ\cr|\cr\circ\cr\ \cr\ \end{matrix}-\circ\dots\circ-
\begin{matrix}\ \circ\cr\ |\cr\ \circ \cr\ \cr\  \end{matrix}-\circ \hskip20mm D_\infty=\bullet-
\begin{matrix}\circ\cr|\cr\circ\cr\ \cr\ \end{matrix}-\circ-\circ\cdots\hskip7mm$$
\vskip-9mm
$$\tilde{E}_6=\bullet-\circ-\begin{matrix}
\circ&\!\!\!\!-\ \circ\cr|&\cr\circ&\!\!\!\!-\ \circ\cr\ \cr\   \end{matrix}-\circ\hskip71mm$$
\vskip-13mm
$$\tilde{E}_7=\bullet-\circ-\circ-
\begin{matrix}\circ\cr|\cr\circ\cr\ \cr\ \end{matrix}-
\circ-\circ-\circ\hskip28mm$$
\vskip-15mm
$$\hskip10mm\tilde{E}_8=\bullet-\circ-\circ-\circ-\circ-
\begin{matrix}\circ\cr|\cr\circ\cr\ \cr\ \end{matrix}-
\circ-\circ$$
\vskip-7mm
$$\hskip62mm A_{\infty}=\bullet-\circ-\circ-\circ\cdots$$
\smallskip

Here the graphs $\tilde{A}_n$ with $n\geq 1$ and $\tilde{D}_n$ with $n\geq 4$ have $n+1$ vertices each. Observe that the graphs $\tilde{E}_n$ with $n=6,7,8$ have $n+1$ vertices each as well.

The first AD graphs are by definition as follows:
$$\tilde{A}_1=\begin{matrix}
\circ\cr
||\cr
\bullet\cr&\cr&\cr\end{matrix}\hskip 2cm \tilde{D}_4=\bullet-\!\!\!\!\!\begin{matrix}
\circ\hskip5mm \circ\cr
\backslash\ \,\slash\cr
\circ\cr&\cr&\cr\end{matrix}\!\!\!\!\!\!\!\!\!\!-\circ$$
\vskip-7mm

The graph $\tilde{A}_n$ as drawn has an even number of vertices. This will be in fact the case: we will only be concerned with graphs of type $\tilde{A}_{2n-1}$, having $2n$ vertices.

We use as well the following graphs, that we call type DE ghost graphs:
$$\tilde{\Delta}_6=\!\!\!\!\begin{matrix}
\circ\hskip5mm \circ\cr
\backslash\ \,\slash\cr
\bullet\cr&\cr&\cr\end{matrix}\!\!\!\!\!\!\!\!\!\!-\circ-
\begin{matrix}\circ\cr|\cr\circ\cr\ \cr\ \end{matrix}-\circ\hskip20mm \tilde{\Delta}_7=
\begin{matrix}\circ\cr|\cr\bullet\cr\ \cr\ \end{matrix}-\circ-
\!\!\begin{matrix}\ \circ\cr\ |\cr\ \circ \cr\ \cr\  \end{matrix}-\circ-\circ-\circ$$
\vskip-9mm

Here the subscript $n=6,7$ comes from the number of vertices, which is $n+1$ as usual. Observe that these two graphs are in fact bad-rooted versions of $\tilde{D}_6$, $\tilde{E}_7$. More comments about them will be given in the end of this section.

The pictures suggest that we have convergences $\tilde{A}_n\to A_{-\infty,\infty}$ and $\tilde{D}_n\to D_\infty$ (and, in some virtual sense, $\tilde{E}_n\to A_\infty$). The following result from \cite{bbi} shows that indeed it is so, at level of circular measures. Moreover, the measures for each series can be recaptured by discretizing the limiting measure at certain roots of unity:

\begin{theorem}
The circular measures of ADE graphs are as follows, where the numbers $t,s$ are given by $t=1,3$ for $n=6,7$ and $s=3,4,6$ for $n=6,7,8$.
\begin{enumerate}
\item For $\tilde{A}_{2n-1}$ we have $\varepsilon=d_n$.
\item For $A_{-\infty,\infty}$ we have $\varepsilon=d$.
\item For $\tilde{D}_{n+2}$ we have $\varepsilon=(d_1'+d_n)/2$.
\item For $D_\infty$ we have $\varepsilon=(d_1'+d)/2$.
\item For $\tilde{\Delta}_n$ we have $\varepsilon=(d_2'+d_t)/2$.
\item For $\tilde{E}_{n}$ we have $\varepsilon=\alpha\,d_s+(d_{s-1}-d_s)/2$.
\item For $A_\infty$ we have $\varepsilon=\alpha\,d$.
\end{enumerate}
\end{theorem}

\begin{proof}
As already mentioned, the result is from \cite{bbi}. We just have to prove the ghost formula (5), not considered in there.

For $\tilde{\Delta}_6$ we denote by $c_k,d_k$ the number of $2k$-paths from the root $c$ to itself, and to the second vertex at right $d$, as indicated below: 
$$\begin{matrix}
\circ\hskip5mm \circ\cr
\backslash\ \,\slash\cr
c\cr&\cr&\cr\end{matrix}\!\!\!\!\!\!\!\!\!\!-\circ-
\begin{matrix}\circ\cr|\cr d\cr\ \cr\ \end{matrix}-\circ$$
\vskip-7mm

Since $c,d$ are at even distances on the graph, the $2k$-paths counted by $c_k,d_k$ can be computed recursively. The equations are as follows:
\begin{eqnarray*}
c_{k+1}&=&3c_k+d_k\cr
d_{k+1}&=&c_k+3d_k
\end{eqnarray*}

The starting values being $c_1=3$, $d_1=1$, we have the following solution:
\begin{eqnarray*}
c_k&=&(4^k+2^k)/2\cr
d_k&=&(4^k-2^k)/2
\end{eqnarray*}

On the other hand, the formula of $c_k$ is the one for multiplicities for the group $D_1$, from previous proof. Thus the algebraic invariants of $\tilde{\Delta}_6$ are those of $D_1$ computed in previous theorem, and in particular we have $\varepsilon=(d_2'+d_1)/2$ as claimed.

For $\tilde{\Delta}_7$ we denote by $c_k,d_k,e_k$ the number of $2k$-paths from the root $c$ to itself, to the second vertex at right $d$, and to the fourth vertex at right $e$: 
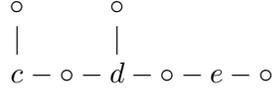
$$\begin{matrix}\circ\cr|\cr c\cr\ \cr\ \end{matrix}-\circ-
\!\!\begin{matrix}\ \circ\cr\ |\cr\ d \cr\ \cr\  \end{matrix}-\circ-e-\circ$$
\vskip-9mm

Since $c,d,e$ are at even distances on the graph, the $2k$-paths counted by $c_k,d_k,e_k$ can be computed recursively. The equations are as follows:
\begin{eqnarray*}
c_{k+1}&=&2c_k+d_k\cr
d_{k+1}&=&c_k+3d_k+e_k\cr
e_{k+1}&=&d_k+2e_k
\end{eqnarray*}

The starting values being $c_1=2$, $d_1=1$, $e_1=0$, we get $c_2=5$, $c_3=15$, $c_4=51$ and so on. The general formula can be deduced by solving the system, and is:
$$c_k=\frac{1}{6}(2+3\cdot 2^k+4^k)$$

On the other hand, the coefficients $m_1=2$, $m_2=3$ and $m=6$ appearing in this formula are nothing but the $m_1,m_2,m$ numbers for the group $S_3$, from previous proof. Thus the algebraic invariants of $\tilde{\Delta}_7$ are those of $S_3$ computed in previous theorem, and in particular we have $\varepsilon=(d_2'+d_3)/2$, and we are done.
\end{proof}

\begin{theorem}
We have the classification table
\setlength{\extrarowheight}{3pt}
\begin{center}
\begin{tabular}[t]{|l|l|l|l|l|l|l|l|l|l|l|}
\hline $\tilde{A}_1$&$\tilde{A}_3$&$\tilde{A}_5$&$\tilde{A}_7$&$\tilde{D}_4$&$\tilde{D}_6$&$\tilde{\Delta}_6$&$\tilde{\Delta}_7$&$\tilde{E}_6$&$\tilde{E}_7$\\ 
\hline $\mathbb Z_1$&$\mathbb Z_2$&$\mathbb Z_3$&$V$&$\mathbb Z_4$&$D_4$&$D_1$&$S_3$&$A_4$&$S_4$\\
\hline
\end{tabular}
\end{center}
\setlength{\extrarowheight}{0pt}
\bigskip
making correspond graphs and groups having the same algebraic invariants.
\end{theorem}

\begin{proof}
This follows from the previous two results.
\end{proof}

We end the section 
with some comments on the following problem: are $\tilde{\Delta}_6$, $\tilde{\Delta}_7$ part of some bigger series of graphs? That is, we would like to have some natural series of graphs $X_n$, having the property that $\tilde{\Delta}_6$, $\tilde{\Delta}_7$ are specializations of it.

There are two possible solutions, coming from the results that we have so far:
\begin{enumerate}
\item $X_n$ is the graph having spectral measure $(e_n+d_n)/2$. 
\item $X_n$ is the graph having spectral measure $(d_2'+d_n)/2$.
\end{enumerate}

The first interpretation, coming from results in previous section, is not the good one: now we know from the general ADE formulae that at even levels we will get usual dihedral graphs, so this series $X_n$ is quite unnatural.

The second interpretation is probably the good one, and we will discuss it here in detail. First, we restrict attention to the ${\rm supp}(\mu)\subset\mathbb N$ case.

\begin{proposition}
There are $5$ measures of type $\varepsilon_n=(d_2'+d_n)/2$ which appear as linear combinations of $\gamma_s$ with $s=0,1,2,3,4$, namely:
\begin{eqnarray*}
\varepsilon_1&=&(\gamma_2+\gamma_4)/2\cr
\varepsilon_2&=&(\gamma_0+2\gamma_2+\gamma_4)/4\cr
\varepsilon_3&=&(2\gamma_1+3\gamma_2+\gamma_4)/6\cr
\varepsilon_4&=&(\gamma_0+6\gamma_2+\gamma_4)/8\cr
\varepsilon_6&=&(\gamma_0+2\gamma_1+6\gamma_2+2\gamma_3+\gamma_4)/12
\end{eqnarray*}
\end{proposition}

\begin{proof}
First we note that
there are $5$ measures of type $d_n$ which appear as linear combinations of $\gamma_s$ with $s=0,1,2,3,4$, namely those corresponding to $n=1,2,3,4,6$,
given by the following formulae:
\begin{eqnarray*}
d_1&=&\gamma_4\cr
d_2&=&(\gamma_0+\gamma_4)/2\cr
d_3&=&(2\gamma_1+\gamma_4)/3\cr
d_4&=&(\gamma_0+2\gamma_2+\gamma_4)/4\cr
d_6&=&(\gamma_0+2\gamma_1+2\gamma_3+\gamma_4)/6
\end{eqnarray*}
The formulae in the proposition are deduced by making averages between $d_2'=\gamma_2$ and the measures $d_n$. The uniqueness assertion is clear.
\end{proof}

We know that the measures $\varepsilon_1,\varepsilon_2,\varepsilon_3$ correspond to the groups $D_1,V,S_3$ and to the graphs $\tilde{\Delta}_6,\tilde{A}_7,\tilde{\Delta}_7$. 

For the remaining measures $\varepsilon_4,\varepsilon_6$, we have a problem
with graphs. For $\varepsilon_4$, the corresponding loop numbers are $c_k=(6\cdot 2^k+4^k)/8$, leading to the following system of equations:
\begin{eqnarray*}
c_{k+1}&=&2c_k+d_k\cr
d_{k+1}&=&c_k+2d_k+3e_k\cr
e_{k+1}&=&d_k+2e_k
\end{eqnarray*}

The starting values here are $c_1=2,d_1=2,e_1=0$, and the first three loop numbers are $c_1=2,c_2=5,c_3=14$. These numeric values, along with the system itself, lead to the conclusion that the graph must look as follows:
$$X_4=\begin{matrix}\circ\cr|\cr \bullet\cr\ \cr\ \end{matrix}-
\circ-\circ-\circ\!\!\!\!\!\begin{matrix}\slash\cr\ \ \  -\!\!\!\!-\!\!\!\!\!\!-\cr\backslash\end{matrix}
\begin{matrix}\!\!\!\!\!\!\!\!\circ-\!\!\!+\!\!\!-\circ\cr\cr
\,\circ-\!\!\!+\!\!\!-\circ\cr\cr\!\!\!\!\!\!\!\!\circ-\!\!\!+\!\!\!-\circ\end{matrix}$$

Here the three special edges at right are negative, in the sense that a loop through any of them adds $-1$ to the number of loops which are counted. The fact that this graph gives the above $c_d,d_k,e_k$ numbers follows by using the method in proof of Theorem 10.1. This is the simplest picture of $X_4$ that we were able to  find. 

The situation for the remaining measure $\varepsilon_6$ is quite similar. As a conclusion, $\tilde{\Delta}_6$, $\tilde{\Delta}_7$ are either exceptional, or are part of a series of non-standard graphs.

\section{Quantum subgroups and classification tables}

We discuss now the non-classical subgroups of $\mathcal Q_4$. Most of them appear as twists of subgroups of $SO(3)$, so we first take a look at these subgroups. These are:
\begin{enumerate}
\item The cyclic groups $\mathbb Z_n$ and the dihedral groups $D_n$.
\item Their limits $SO(2)$ and $O(2)$, plus $SO(3)$ itself.
\item The exceptional groups $A_4,S_4,A_5$. 
\end{enumerate}

The McKay correspondence \cite{mc} provides us with an ADE classification of these subgroups. For the purposes of this paper, we will need a direct proof of this result, in terms of algebraic invariants similar to those considered in previous sections.

We use the canonical identification $SO(3)\simeq PU(2)$, which endows $SO(3)$ with a fundamental representation $\pi$, on the Hilbert space $M_2(\mathbb C)\simeq \mathbb C^4$. We have $\pi=1+\rho$, where $\rho$ is the canonical representation on $\mathbb C^3$, coming from that on $\mathbb R^3$.

\begin{definition}
The algebraic invariants of a subgroup $G\subset SO(3)$ are:
\begin{enumerate}
\item The multiplicity $c_k$. This is the number of copies of $1$ into $\pi^{\otimes k}$.
\item The spectral measure $\mu$. This is the measure having $c_k$ as moments.
\item The circular measure $\varepsilon$. This is the pullback of $\mu$ via $\varphi(q)=(q+\bar{q})^2$.
\end{enumerate}
\end{definition}

In this definition $\mu$ and $\varepsilon$ are as usual probability measures on the real line, and on the unit circle.

With these notations, the ADE classification of subgroups of $SO(3)$ coming from the McKay correspondence can be stated in the following way:

\begin{theorem}
We have the classification table
\setlength{\extrarowheight}{3pt}
\begin{center}
\begin{tabular}[t]{|l|l|l|l|l|l|l|l|l|l|l|}
\hline $\tilde{A}_{2n-1}$&$A_{-\infty,\infty}$&$A_{\infty}$&$\tilde{D}_{n+2}$&$D_\infty$&$\tilde{E}_6$&$\tilde{E}_7$&$\tilde{E}_8$\\ 
\hline $\mathbb Z_n$&$SO(2)$&$SO(3)$&$D_{n}$&$O(2)$&$A_4$&$S_4$&$A_5$\\ 
\hline
\end{tabular}
\end{center}
\setlength{\extrarowheight}{0pt}
\bigskip
making correspond graphs and groups having the same algebraic invariants.
\end{theorem}

\begin{proof}
It is convenient to use subgroups of $SU(2)$ rather than subgroups of $SO(3)$. Recall first that we have the following equality:
$$SU(2)=\left\{ \begin{pmatrix}a&b\cr -\bar{b}&\bar{a}\end{pmatrix}\,\Big\vert\, |a|^2+|b|^2=1\right\}$$

We regard $a,b$ as complex functions on $SU(2)$, or as elements of ${\rm C}(SU(2))$. We also use the real functions $x,y,z,t$ given by 
$a=x+iy$ and  $b=z+it$.

The change of variables $(a,b)\to (x,y,z,t)$ gives an isomorphism of real algebraic varieties $SU(2)\simeq S^3$. This isomorphism is measure-preserving, in the sense that the Haar measure on $SU(2)$ corresponds to the uniform measure on $S^3$.

Now given a subgroup $G\subset SO(3)$, we take its double cover $G'\subset SU(2)$. The multiplicity $c_k$ for the group $G$, as constructed in the above definition, is
$$c_k=\int\chi^{2k}
=\int(a+\bar{a})^{2k}
=\int (4x^2)^{k}$$

Here integration is with respect to the Haar measure of $G'$, and we use the functions $\chi,a,x\in {\rm C}(G')$, where $\chi$ is the character of the fundamental representation. 

We can compute now the invariants in the statement.

For $G=SO(3)$ we have $G'=SU(2)$. Since $x$ is semicircular on $[-1,1]$, its double $2x$ is semicircular on $[-2,2]$, so the quadrupled square $4x^2$ is free Poisson on $[0,4]$. We get the measure $\varepsilon=\alpha\,d$, corresponding to the graph $A_\infty$.

For $G=\mathbb Z_n$ we have $G'=\mathbb Z_{2n}$, given in matrix form by:
$$\mathbb Z_{2n}=\left\{ \begin{pmatrix}q&0\cr 0&\bar{q}\end{pmatrix}\,\Big\vert\, q^{2n}=1\right\}$$

This gives the following formula, where the sum is over $2n$-roots of unity:
$$c_k=\frac{1}{2n}\sum_q(q+\bar{q})^{2k}$$

We get the circular measure $\varepsilon=d_n$, corresponding to the graph $\tilde{A}_{2n-1}$.

For $G=D_n$ we have $G'=D_{2n}$, given in matrix form by:
$$D_{2n}=\left\{ \begin{pmatrix}q&0\cr 0&\bar{q}\end{pmatrix}, \begin{pmatrix}0&q\cr -\bar{q}&0\end{pmatrix}\,\Big\vert\, q^{2n}=1\right\}$$

This gives the following formula, where the sum is over $2n$-roots of unity:
$$c_k=\frac{1}{4n}\sum_q(q+\bar{q})^{2k}$$

We get the circular measure $\varepsilon=(d_1'+d_n)/2$, corresponding to the graph $\tilde{D}_{n+2}$.

The limiting cases $SO(2)$ and $O(2)$ are obtained by taking limits.

For $A_4,S_4,A_5$ the situation is a bit more complicated. In lack of a simple direct proof, we prefer to refer here to the McKay correspondence.
\end{proof}

Our ADE table is different from the usual one, because it is based on projective representations. For instance the cyclic groups $\mathbb Z_n$, usually labeled $\tilde{A}_n$, are here subject to a rescaling, make them correspond to $\tilde{A}_{2n-1}$. This is not very surprising, say because of the subfactor point of view, very close to ours: $\tilde{A}_{2n}$ is known not to be a principal graph. See \cite{ek}.

\medskip

We are now in position of finishing the ADE classification of subgroups of
$\mathcal Q_4$. Recall from Theorem 8.1 and Theorem 10.2 that there are a
number of subgroups left: 3 infinite series and 5 single quantum subgroups.

\begin{theorem}
We have the classification table
\setlength{\extrarowheight}{3pt}
\begin{center}
\begin{tabular}[t]{|l|l|l|l|l|l|l|l|l|l|l|}
\hline $\tilde{A}_{4n-1}$&$A_{-\infty,\infty}$&$A_{\infty}$&$\tilde{D}_{2n+2}$&$D_\infty$&$\tilde{E}_7$&$\tilde{E}_8$\\ 
\hline $\widehat{D}_n$& $\widehat{D}_\infty$&$\mathcal Q_4$&$D_{2n}^\tau$, $DC_n^\tau$ &$O_{-1}(2)$&$S_4^\tau$&$A_5^\tau$\\
\hline
\end{tabular}
\end{center}
\setlength{\extrarowheight}{0pt}
\medskip
making correspond graphs and quantum groups having the same algebraic invariants.
\end{theorem}

\begin{proof}
It is clear that the algebraic invariants are representation theoretic invariants,
and remain the same when passing from a subgroup of $SO(3)$ to its twisted (or pseudo-twisted)
version in $SO(3)$. Therefore, for the twists $\Q_4 \cong SO_{-1}(3)$,
$O_{-1}(2)$, $D_{2n}^\tau$, $DC_n^\tau$, $S_4^\tau$, $A_5^\tau$
we can use the previous ADE classification of subgroups of $SO(3)$.

For $\widehat{D}_n$ with $n=3,4,\ldots,\infty$ we have to compute the
invariants of a certain $4$-dimensional representation. But, as for any group
dual, this representation must be of the form $\pi={\rm diag}(g_i)$, with
$D_n=<g_1,g_2,g_3,g_4>$. Moreover, the multiplicities are given by
$c_k=\#\{(i_1,\ldots,i_k)\mid g_{i_1}\ldots g_{i_k}=1\}$. Here we have 
$D_n = \langle  1, g, 1, h \rangle$ where $g^2=1=h^2$ are the given
generators of $D_n$ (with the additionnal relation $(gh)^n=1$ if $n$ is finite). Thus we
have to compute the $k$-loops in the Cayley graph of $D_n$ with respect to the
set $\{1,g,1,h\}$, which is $\tilde{A}_{2n-1}$ with two $1$-loops added at
each vertex. But this is exactly the number of $2k$-loops of
$\tilde{A}_{4n-1}$, and this gives the result.    
\end{proof}


The results in previous sections give the following classification table, where
the relation between graphs and subgroups is that the algebraic invariants are
the same. The last column contains the idea of the proof.

{\small\setlength{\extrarowheight}{3pt}
\begin{center}
\begin{tabular}[t]{|l|l|l|l|l|l|}
\hline ADE graph&Subgroup of $\mathcal Q_4$&Correspondence\\
\hline $\tilde{A}_1$, $\tilde{A}_3$, $\tilde{A}_5$, $\tilde{A}_7$&$\mathbb Z_1$, $\mathbb Z_2$, $\mathbb Z_3$, $V$&ADE for $S_4$\\
\hline $\tilde{A}_{4n-1}$ & $\widehat{D}_n$ & Computation\\
\hline $A_{-\infty,\infty}$&$\widehat{D}_\infty$&Computation\\
\hline $A_\infty$&$\mathcal Q_4$&Twisted McKay\\
\hline $\tilde{D}_4$, $\tilde{D}_6$&$\mathbb Z_4$, $D_4$&ADE for $S_4$\\
\hline $\tilde{D}_{2n+2}$&$D_{2n}^\tau$, $DC_n^\tau$&Twisted McKay\\
\hline $D_\infty$&$O_{-1}(2)$&Twisted McKay\\
\hline $\tilde{\Delta}_6$, $\tilde{\Delta}_7$&$D_1$, $S_3$&ADE for $S_4$\\
\hline $\tilde{E}_6$&$A_4$&ADE for $S_4$\\
\hline $\tilde{E}_7$&$S_4$, $S_4^\tau$&McKay + twist\\
\hline $\tilde{E}_8$&$A_5^\tau$&Twisted McKay\\
\hline
\end{tabular}
\end{center}
\setlength{\extrarowheight}{0pt}}

\medskip

A second classification table can be obtained by collecting together the ADE
classification information for $\mathcal Q_4$ and for $SO(3)$ (the squares
denote missing subgroups):
{\small\setlength{\extrarowheight}{3pt}
\begin{center}
\begin{tabular}[t]{|l|l|l|l|l|l|}
\hline Subgroup of $SO(3)$&ADE graph&Subgroup of $\mathcal Q_4$\\
\hline $\mathbb Z_1$, $\mathbb Z_2$, $\mathbb Z_3$, $\mathbb Z_4$&$\tilde{A}_1$, $\tilde{A}_3$, $\tilde{A}_5$, $\tilde{A}_7$&$\mathbb Z_1$, $\mathbb Z_2$, $\mathbb Z_3$, $V$\\
\hline $\mathbb Z_{2n-1}$, $\mathbb Z_{2n}$&$\tilde{A}_{4n-3}$,
$\tilde{A}_{4n-1}$&$\square$, $\widehat{D}_n$\\
\hline $SO(2)$&$A_{-\infty,\infty}$&$\widehat{D}_\infty$\\
\hline $SO(3)$&$A_\infty$&$\mathcal Q_4$\\
\hline $V$, $D_4$&$\tilde{D}_4$, $\tilde{D}_6$&$\mathbb Z_4$, $D_4$\\
\hline $D_{2n-1}$, $D_{2n}$&$\tilde{D}_{2n+1}$, $\tilde{D}_{2n+2}$& $\square$, $D_{2n}^\tau$, $DC_n^\tau$\\
\hline $O(2)$&$D_\infty$&$O_{-1}(2)$\\
\hline $\square$, $\square$ &$\tilde{\Delta}_6$, $\tilde{\Delta}_7$&$D_1$, $S_3$\\
\hline $A_4$&$\tilde{E}_6$&$A_4$\\
\hline $S_4$&$\tilde{E}_7$&$S_4$, $S_4^\tau$\\
\hline $A_5$&$\tilde{E}_8$&$A_5^\tau$\\
\hline
\end{tabular}
\end{center}
\setlength{\extrarowheight}{0pt}}
\bigskip

The following problem appears: is there any ``spatial'' correspondence between ADE graphs and subgroups of $\mathcal Q_4$, beyond representation theory results? 

There are several interpretations here:
\begin{enumerate}
\item McKay correspondence. There are many ways of establishing the McKay correspondence for $SO(3)$, see \cite{mc}. However, none of them seems to have a clear extension to $\mathcal Q_4$, or at least we don't know how to find it. One problem here is to find the correct quantum analogue of Platonic solids.
\item Planar algebras. A statement a bit more ``spatial'' than the representation theory one can be formulated in terms of planar algebras, by using results in \cite{ba2}. Thus our ADE classification can be regarded as being part of the ADE classification of index $4$ planar algebras, discussed in \cite{jo2}.
\item Other. There are many other ADE results, for instance Arnold's celebrated classification of singularities \cite{ar}. Finding a possible relation with quantum groups seems to be a key problem, that we would like to raise here.
\end{enumerate}

\section{Concluding remarks}

After Wang's discovery of the quantum permutation group
$\Q_n$, the classification of its subgroups at $n=4$ given here
strengthens the idea that 
quite unexpected quantum groups might act on classical spaces. 
In view of Klein's program, this means that a very simple classical space, such as the space
consisting of $4$ points, might have a very rich quantum geometry.

\smallskip

A natural fundamental question here concerns the existence and construction
of Woronowicz' differential calculi \cite{wodc} on the quantum permutation groups listed in Theorem 1.1.
With the exception of the group algebra $\mathbb C[D_\infty]$, all 
the corresponding Hopf algebras are coquasitriangular, so there are well known
methods for constructing differential calculi on them. 
We have two different families:

$\bullet$ The Hopf algebras $\r(\Q_4)$, $\r(O_{-1}(2))$, $\r(G)$ for $G \subset S_4$,
$\r(S_4^\tau)$, $\r(A_5^\tau)$, $\r(D_n^\tau)$. These Hopf algebras
are cotriangular because they are cocycle deformations of function
algebras on classical groups. A result of Majid and Oeckl \cite{mo}
ensures that there is a bijective correspondence between bicovariant differential calculi on these Hopf algebras
and on the corresponding function algebras.

$\bullet$ The Hopf algebras $\r(DC_n^\tau)$. These are not cotriangular 
but are coquasitriangular by Suzuki's paper \cite{su}. Therefore we can use the general method
of Klimyk and Schm\"udgen
in chapter 14 of the book \cite{ks} for constructing a differential
calculus on them.

\smallskip

Concerning the classification problem for quantum permutation groups, 
it is probably only reachable for small 
$n$, say $n=5,6,7$, and for general $n$ we only hope to 
have results splitting the family
of quantum subgroups of $\Q_n$ in subfamilies having
well identified properties.
As already said in the introduction, a first 
obstruction is the complexity of the $\c^*$-algebra 
$\c(\Q_n)$, because of the non-amenability
of the discrete quantum group dual to $\Q_n$ \cite{ba0}.

A first easy thing to be done is to generalize the usual notions of
multiple transitivity to quantum permutation groups, using
algebraic invariants as in Section 8. 
We would like to mention the following question on multiple
transitivity. It is known that a classical 
$k$-transitive permutation group with $k\geq 6$
is either a symmetric or an alternating group. 
Is there such a kind of result for quantum permutation groups?

Also it will also be interesting to examine
the quantum permutation groups on a prime number of points $p$. 
Recall that a classical result of Burnside asserts that if 
$G$ is a subgroup of $S_p$ acting transitively on $p$ points, then
$G$ is either 2-transitive or is a proper subgroup of the affine group
$AGL_1(\mathbb Z_p)$ \cite{dm}. 
Is there an analogue of Burnside's theorem for quantum permutation groups?
A related question was studied in \cite{bbch}.

Finally we have the following conjecture.

\begin{conjecture}
The  symmetric group $S_n$ is maximal in the quantum permutation group $\Q_n$,
for any $n$. 
\end{conjecture}

As a direct consequence of the classification theorem of quantum subgroups of 
$\Q_4$, the conjecture is true at $n=4$.


\begin{thebibliography}{99}

\bibitem{ar} V.I. Arnold, S.M. Gusein-Zade and A.N. Varchenko, Singularities of differentiable maps, Birkh\"auser (1985).

\bibitem{ba0} T. Banica, Symmetries of a generic coaction,
{\em Math. Ann.} {\bf 314} (1999), 763-780.

\bibitem{ba1} T. Banica, Subfactors associated to compact Kac algebras,
{\em Integral Equations Operator Theory} {\bf 39} (2001), 1-14. 

\bibitem{ba2} T. Banica, Quantum automorphism groups of homogeneous graphs, {\em J. Funct. Anal.} {\bf 224} (2005), 243-280.

\bibitem{bb} T. Banica and J. Bichon, Free product formulae for quantum
  permutation groups, {\em J. Inst. Math. Jussieu} {\bf 6} (2007), 381-414.
.

\bibitem{bb2} T. Banica and J. Bichon, Quantum automorphism groups of
  vertex-transitive graphs of order $\leq 11$, {\em J. Algebraic Combin.} {\bf
    26} (2007), 83-105.

\bibitem{bbch} T. Banica, J. Bichon and G. Chenevier, Graphs having no quantum
  symmetry, {\em Ann. Inst. Fourier} {\bf 57} (2007), 955-971.

\bibitem{bbc} T. Banica, J. Bichon and B. Collins, The hyperoctahedral quantum group, {\tt math.RT/0701859}.

\bibitem{bbi} T. Banica and D. Bisch, Spectral measures of small index principal graphs, {\em Comm. Math. Phys.} {\bf 269} (2007), 259-281.

\bibitem{bc} T. Banica and B. Collins, Integration over quantum permutation groups, {\em J. Funct. Anal.} {\bf 242} (2007), 641-657.

\bibitem{bc2} T. Banica and B. Collins, Integration over the Pauli quantum group,
\texttt{math.QA/0610041}.

\bibitem{bm} T. Banica and S. Moroianu, On the structure of quantum permutation groups, {\em Proc. Amer. Math. Soc.} {\bf 135} (2007), 21-29. 

\bibitem{bi0} J. Bichon, Quelques nouvelles d\'eformations du groupe sym\'etrique, {\em C. R. Acad. Sci. Paris} {\bf 330} (2000), 761-764.

\bibitem{bi1} J. Bichon, Quantum automorphism groups of finite graphs, {\em Proc. Amer. Math. Soc.} {\bf 131} (2003), 665-673.

\bibitem{bi2} J. Bichon, Free wreath product by the quantum permutation group, {\em Alg. Rep. Theory} {\bf 7} (2004), 343-362.

\bibitem{da} A. Davydov, Galois algebras and monoidal functors between categories of representations of finite groups, {\em J. Algebra} {\bf 244} (2001), 273-301.

\bibitem{dm} J. Dixon and B. Mortimer, Permutation groups, GTM 163, Springer-Verlag (1996).

\bibitem{do} Y. Doi, Braided bialgebras and quadratic algebras, {\em Comm. Algebra} {\bf 21} (1993), 1731-1749.

\bibitem{dr} V. Drinfeld, Quantum groups, Proc. ICM Berkeley (1986), 798-820.

\bibitem{ev}M. Enock and L. Vainerman, Deformation of a Kac algebra by an abelian subgroup, {\em Comm. Math. Phys.} {\bf 178} (1996), 571--596.

\bibitem{eg1} P. Etingof and S. Gelaki, Some properties of finite-dimensional semisimple Hopf algebras, {\em Math. Res. Lett.}  {\bf 5} (1998), 191-197.

\bibitem{eg2} P. Etingof and S. Gelaki, The representation theory of cotriangular semisimple Hopf algebras, {\em Internat. Math. Res. Notices}  {\bf 7} (1999), 387-394.

\bibitem{eg3} P. Etingof and S. Gelaki, The classification of triangular semisimple and cosemisimple Hopf algebras over an algebraically closed field, {\em Internat. Math. Res. Notices}  {\bf 5} (2000), 223-234.

\bibitem{ek} D. Evans and Y. Kawahigashi, Quantum symmetries on operator algebras, Oxford University Press (1998).


\bibitem{jo2} V.F.R. Jones, The annular structure of subfactors, {\em Monogr. Enseign. Math.} {\bf 38} (2001), 401-463.

\bibitem{ks} A. Klimyk and K. Schm\"udgen, Quantum groups and their representations, Texts and Monographs in Physics, Springer-Verlag, Berlin (1997).

\bibitem{mo} S. Majid and R. Oeckl, Twisting of quantum differentials and the
  Planck scale Hopf algebra,  {\em Comm. Math. Phys.}  {\bf 205}  (1999), 617-655. 

\bibitem{ma} Y. Manin, Quantum groups and noncommutative geometry, {\em Publications du CRM} {\bf 1561}, Univ. de Montr\'eal (1988).

\bibitem{mp} V.A. Marchenko and L.A. Pastur, Distribution of eigenvalues in certain sets of random matrices, {\em Mat. Sb.} {\bf 72} (1967), 507-536.

\bibitem{mas} A. Masuoka, Cocycle deformations and Galois objects for some cosemisimple Hopf algebras of finite dimension, {\em Contemp. Math.} {\bf 267} (2000), 195-214. 

\bibitem{mc} J. McKay, Graphs, singularities and finite groups, {\em Proc. Symp. Pure Math.} {\bf 37} (1980), 183-186.

\bibitem{nz} W.D. Nichols and  M.B. Zoeller, A Hopf algebra freeness theorem, {\em Amer. J. Math.} {\bf 111}  (1989), 381-385. 

\bibitem{ni} D. Nikshych, $K_0$-rings and twisting of finite-dimensional semisimple Hopf algebras, {\em Comm. Algebra} {\bf 26} (1998), 321-342.

\bibitem{po} P. Podles, Symmetries of quantum spaces. Subgroups and quotient spaces of quantum ${\rm SU}(2)$ and ${\rm SO}(3)$ groups, {\em Comm. Math. Phys.} {\bf 170} (1995), 1-20. 

\bibitem{ro} M. Rosso, Alg\`ebres enveloppantes quantifi\'ees, groupes quantiques compacts de matrices et calcul differentiel non-commutatif, {\em Duke Math. J.} {\bf 61} (1990), 11-40.

\bibitem{sc1} P. Schauenburg, Hopf bigalois extensions, {\em Comm. Algebra} {\bf 24} (1996), 3797-3825.

\bibitem{su} S. Suzuki, A family of braided cosemisimple Hopf algebras of finite dimension,
{\em Tsukuba J. Math.} {\bf 22} (1998), 1-29. 

\bibitem{va} L. Vainerman, $2$-cocycles and twisting of Kac algebras,   {\em Comm. Math. Phys.} {\bf 191} (1998), 697-721.

\bibitem{vdn} D.V. Voiculescu, K.J. Dykema and A. Nica, Free random variables, {\em CRM Monograph Series} {\bf 1}, AMS (1992).

\bibitem{wa0} S. Wang, Free products of compact quantum groups, {\em Comm. Math. Phys.} {\bf 167} (1995), 671-692.

\bibitem{wa} S. Wang, Quantum symmetry groups of finite spaces, {\em Comm. Math. Phys.} {\bf 195} (1998), 195-211.

\bibitem{wo} S.L. Woronowicz, Compact matrix pseudogroups, {\em Comm. Math. Phys.} {\bf 111} (1987), 613-665.

\bibitem{wodc} S.L. Woronowicz, Differential calculus on compact matrix
  pseudogroups (quantum groups),  {\em Comm. Math. Phys.}  {\bf 122} (1989), 125-170.

\bibitem{wo2} S.L. Woronowicz, Compact quantum groups, in  {\em ``Sym\'etries quantiques'' (Les Houches, 1995)}, North Holland, Amsterdam (1998), 845-884.

\bibitem{z} S. Zakrzewski, Matrix pseudogroups associated with anti-commutative plane, {\em Lett. Math. Phys.} {\bf 21} (1991), 309-321. 

\end{thebibliography}
\end{document}